\documentclass[11pt]{article}%
\usepackage{amssymb}
\usepackage{amsmath}
\usepackage{amsfonts}
\usepackage{boxedminipage}
\usepackage{color}
\usepackage{fancyhdr}
\usepackage{makeidx}
\usepackage{graphicx}
\usepackage{geometry}
\usepackage{multicol}%
\providecommand{\U}[1]{\protect\rule{.1in}{.1in}}
\newtheorem{proposition}{Proposition}[section]
\newtheorem{theorem}[proposition]{Theorem}

\newtheorem{definition}[proposition]{Definition}
\newtheorem{remark}[proposition]{Remark}
\newtheorem{example}[proposition]{Example}

\numberwithin{equation}{section}
\numberwithin{proposition}{section}

\makeindex
\begin{document}

\title{Importance sampling for metastable and multiscale dynamical systems}
\author{K. Spiliopoulos\footnote{Department of  Mathematics and Statistics,
Boston University, Boston, MA, 02215, kspiliop@math.bu.edu. This work was partially supported by the National Science Foundation
CAREER award DMS 1550918 }}
\date{\today}

\maketitle

\begin{abstract}
In this article, we address the issues that come up in the design of importance sampling schemes for rare events associated to stochastic dynamical systems. We focus on the issue of metastability and on the effect of multiple scales.  We discuss why seemingly reasonable schemes that follow large deviations optimal paths may perform poorly in practice, even though they are asymptotically optimal. Pre-asymptotic optimality is important when one deals with metastable dynamics  and we discuss possible ways as to how to address this issue. Moreover, we discuss how the effect of the multiple scales (either in periodic or random environments) on the efficient design of importance sampling should be addressed. We discuss the mathematical and practical issues that come up, how to overcome some of the issues and discuss future challenges.
\end{abstract}

\section{Introduction}\label{S:Introduction}

%\section{Rare events for small noise diffusions and metastability}\{]

In this paper, we discuss recent developments on importance sampling methods for metastable dynamics that may also have multiple scales. Development of accelerated Monte Carlo methods for metastable, multiple-scale processes is of great interest.  Importance sampling is a variance reduction technique in Monte-Carlo simulation, which is especially relevant when dealing with rare events. Since its introduction, importance sampling has been one of the most popular techniques for rare event simulation. There is a vast literature of papers investigating its applications from a broad set of sciences including engineering, chemistry, physics, biology, finance, insurance, e.g., \cite{AsmussenGlynn2007,BoyleBroadieGlasserman,Glasserman2004,GlynnIglehart1988,GriffithsTavare1994,Levine,MazonkaJarzynskiBlocki1998,Siegmund1976,VielPatelNiyazWhaley2002,ZuckermanWoolf}.

Consider a sequence $\{X^{\epsilon}\}_{\epsilon>0}$ of random elements and assume that we want to estimate the probability $0<p^{\epsilon}=\mathbb{P}\left[X^{\epsilon}\notin \mathcal{D}\cup \partial \mathcal{D}\right]\ll 1$ for a given set $\mathcal{D}$, such that the event $\left\{X^{\epsilon}\notin \mathcal{D}\cup \partial \mathcal{D}\right\}$ is unlikely for small $\epsilon$. If closed form formulas are not available, or numerical approximations are either too crude or unavailable, then one has to resort in simulation. It is well known that standard Monte-Carlo simulation techniques (i.e., using the unbiased estimator $\hat{p}^{\epsilon}=\frac{1}{N}\sum_{j=1}^{N}1_{X^{\epsilon,j}\notin \mathcal{D}\cup \partial \mathcal{D}}$) perform rather poorly in the rare-event regime. As it is known, see for example \cite{AsmussenGlynn2007},  in order to  achieve relative error smaller than one using standard Monte Carlo, one needs an effective sample size $N\approx 1/p^{\epsilon}$. In other words, for a fixed computational
cost, relative errors grow rapidly as the event becomes more rare.   Thus  standard Monte-Carlo is infeasible for rare-event simulation.

The goal of importance sampling is to simulate the system under an
alternative probability distribution $\bar{\mathbb{P}}$ instead of the original probability $\mathbb{P}$. Let's say for example that we are interested in the estimation of
 \begin{equation}
\mathbb{E}_{y}[e^{-\frac{1}{\epsilon}h(X^{\epsilon}_{T})}]\textrm{
or } \mathbb{P}_{y}\left[\tau^{\epsilon}_{\mathcal{D}\cup \partial \mathcal{D}}\leq T\right]\label{Eq:ProbToEstimate}
\end{equation}
where $h:\mathbb{R}^{d}\mapsto\mathbb{R}$ is a positive function, $T>0, \epsilon>0$, $y\in D$ is the initial point, $\tau^{\epsilon}_{\mathcal{D}\cup \partial \mathcal{D}}$ is exit time from the set $\mathcal{D}\cup \partial \mathcal{D}$, $X^{\epsilon}$ is a stochastic process modeling the dynamics.  Also, notice that the probability  above can be considered (modulo the important technical point of lack of continuity) as a special case of $\mathbb{E}_{y}[e^{-\frac{1}{\epsilon}h(X^{\epsilon}_{T})}]$, when $h$ is for example chosen such that $h(x)=0$ for $x\notin \mathcal{D}\cup \partial \mathcal{D}$ and $h(x)=+\infty$ for $x\in \mathcal{D}\cup \partial \mathcal{D}$.

When rare events dominate, then standard Monte-Carlo
methods perform poorly in the small noise limit. %%Apart from the
%%smallness of the noise, an additional difficulty  is the presence
%%of the fast
%% oscillating coefficients.
% Importance sampling is a variance reduction technique in Monte-Carlo simulation. The goal is to simulate the system under an
%alternative probability distribution $\bar{\mathbb{P}}$
%instead of the original probability $\mathbb{P}$.
 Then, to estimate  $\mathbb{E}_{y}[e^{-\frac{1}{\epsilon}h(X^{\epsilon}_{T})}]$, one generates iid samples $X^{\epsilon}_{(k)}$ from $\bar{\mathbb{P}}$ and
uses the importance sampling estimator %$\frac{1}{K} \sum_{k=1}^K  e^{-\frac{1}{\epsilon}h(X_k)} (d\mathbb{P}/d\bar{\mathbb{P}})(X_k)$.
\begin{equation}
\frac{1}{N} \sum_{k=1}^N  e^{-\frac{1}{\epsilon}h(X^{\epsilon}_{(k)})} \frac{d\mathbb{P}}{d\bar{\mathbb{P}}}(X^{\epsilon}_{(k)}).\label{Eq:MonteCarloEstimate}
\end{equation}
The key question is the design of $\bar{\mathbb{P}}$ such that the second moment $\bar{\mathbb{E}}_{y}[e^{-\frac{1}{\epsilon}h(X^{\epsilon}_{T})}(d\mathbb{P}/d\bar{\mathbb{P}})(X^{\epsilon}_{\cdot})]^{2}$
% \begin{equation}
% \bar{\mathbb{E}}_{y}\left[e^{-\frac{1}{\epsilon}h(X^{\epsilon}_{T})}
% \frac{d\mathbb{P}}{d\bar{\mathbb{P}}}(X^{\epsilon}_{\cdot})\right]^{2}\label{Eq:MinimumVariance}
% \end{equation}
(and hence the variance) is minimized. $\bar{\mathbb{E}}$ is
the expectation operator under $\bar{\mathbb{P}}$. The choice of the appropriate alternative measure $\bar{\mathbb{P}}$  is closely related to certain Hamilton-Jacobi-Bellman (HJB) equations.

%The purpose of this article is to summarize recent developments and point out issues and challenges in the design of provably efficient importance sampling schemes for stochastic dynamical systems that may have metastable characteristics and multiple scales.
The first issue that we address is the effect of rest points (and metastability in general) on importance sampling. It turns out that when dealing with metastability, even seemingly reasonable schemes that are also asymptotically optimal, may perform poorly in practice. This includes also changes of measure that try to enforce the simulated trajectories to follow large deviations most likely paths. The reason for the degradation in performance is the role of prefactors. Prefactors can become very important when rest points are included in the domain of interest for the simulation. Large deviations based change of measures may not account for the prefactors, as they rely on logarithmic asymptotics. We elaborate on these issues and discuss potential ways on how the issue can be addressed.

The second issue that we address is the effect of multiple scales on the design of provably-efficient importance sampling methods. It turns out that when the dynamical system has widely separated multiple scales, then one can use averaging and homogenization techniques. However, as we shall see, it is not sufficient to base the design of importance sampling on the effective homogenized dynamics. The local information needs to be taken into account. Mathematically this is done using the so called cell problem, or macroscopic problem, in the theory of periodic and random homogenization.

The rest of the article is summarized as follows. In Section \ref{S:ReviewIS_LDP} we review the classical large deviations theory and the setup of importance sampling for small noise diffusions. In Section \ref{S:EffectRestPoints} we discuss the effects of rest points, i.e. of stable and unstable equilibrium points, in the design of importance sampling. We argue why asymptotic optimality may actually not mean good practical performance and we also argue that following large deviations most likely optimal paths may lead to poor performance. In addition, we present constructions that lead to guaranteed good performance. We supplement the theoretical arguments by simulation studies. We refer the interested reader to \cite{DupuisSpiliopoulosZhou2013,DupuisSpiliopoulos2014a} for more details. In Sections \ref{S:IS_MultiscaleLangevin} and \ref{S:RareEventRandomEnvironments}, we address the design of importance sampling schemes in the presence of multiple scales. We construct asymptotically optimal schemes in the presence of multiple scales. To be more precise, in Section \ref{S:IS_MultiscaleLangevin} we consider overdamped Langevin dynamics in periodic multiscale environments and we review the related large deviations theory and importance sampling theory, presenting simulation studies.  The interested reader can also consult \cite{DupuisSpiliopoulosWang, DupuisSpiliopoulosWang2}. In Section \ref{S:RareEventRandomEnvironments} we review recent developments in large deviations and importance sampling for multiscale dynamics in random environments, see also \cite{Spiliopoulos2013c,Spiliopoulos2015}. In Section \ref{S:IS_metastabilityMultiscale} we describe how one can combine the results of Section
\ref{S:EffectRestPoints} with those of Sections \ref{S:IS_MultiscaleLangevin} and \ref{S:RareEventRandomEnvironments} and also review future directions.

For the sake of concreteness and for exposition purposes we restrict the presentation of this article in the case of diffusions with gradient drift and constant diffusivity, which also implies reversible diffusion dynamics. However, we mention that almost all of the arguments can and have been generalized to the case with general state dependent drift and diffusion coefficient, especially those about the effect of multiple scales on importance sampling, see \cite{DupuisSpiliopoulos,DupuisSpiliopoulosWang, DupuisSpiliopoulosZhou2013,Spiliopoulos2013c,Spiliopoulos2015}. For results in the infinitely dimensional case we refer the interested reader to \cite{SalinsSpiliopoulos2016}.

\section{Review of large deviations and importance sampling theory for diffusions}\label{S:ReviewIS_LDP}
Let us briefly review the setup for small noise diffusions in $\mathbb{R}^{d}$ (e.g. \cite{DupuisSpiliopoulosWang,VandenEijndenWeare}) \textit{without} the effect of multiple scales. Let $W_{t}$ be a standard d-dimensional Wiener process and consider
\begin{equation}
dX_{t}^{\epsilon }=-\nabla V( X_{t}^{\epsilon })  dt+\sqrt{\epsilon }\Gamma dW_{t},\hspace{%
0.2cm}X_{t_{0}}^{\epsilon }=y.\label{Eq:LDPnoHomogenization}
\end{equation}

%Let us define $a(x)=\sigma(x)\sigma^{T}(x)$.
Large deviations principle for the process $X^{\epsilon}_{t}$ is well known (e.g, \cite{FreidlinWentzell88}). In particular, the action functional  for the process
$X^{\epsilon}_{t}, t_{0}\leq t\leq T$, in $\mathcal{C}([t_{0},T])$ as $\epsilon
\downarrow 0$ has the form $\frac{1}{\epsilon}S_{t_{0}T}(\phi)$, where %$ S_{t_{0}T}(\phi)=+\infty$ for $\phi\notin \mathcal{AC}([t_{0},T])$ and
%\begin{equation}
% S_{t_{0}T}(\phi)=\frac{1}{2}\int_{t_{0}}^{T}(\dot{\phi}_{s}+\nabla V(\phi_{s}))^{T}(\Gamma\Gamma^{T})^{-1}(\dot{\phi}_{s}+\nabla V(\phi_{s}))ds,  \text{ if } \phi\in \mathcal{AC}([t_{0},T])
%\label{ActionFunctionalSmooth}
%\end{equation}

\begin{eqnarray}
 S_{t_{0}T}(\phi)=\begin{cases}\frac{1}{2}\int_{t_{0}}^{T}(\dot{\phi}_{s}+\nabla V(\phi_{s}))^{T}\left[\Gamma\Gamma^{T}\right]^{-1}(\dot{\phi}_{s}+\nabla V(\phi_{s}))ds, & \text{if } \phi\in \mathcal{AC}([t_{0},T])
 \\
               +\infty, & \text{otherwise.}
       \end{cases}  \label{ActionFunctionalSmooth}
\end{eqnarray}
Here $\mathcal{C}([t_{0},T])$, $\mathcal{AC}([t_{0},T])$ are the sets of continuous and absolutely continuous functions on $[t_{0},T]$ respectively.
%This allows to characterize events such as $\mathbb{P}_{y}\left[X^{\epsilon}_{T}\in A\right]$.
Then, under fairly general conditions,
\begin{equation*}
\mathbb{E}_{y}[e^{-\frac{1}{\epsilon}h(X^{\epsilon}_{T})}]\approx e^{-\frac{1}{\epsilon}\inf\left\{S_{t_{0}T}(\phi)+h(\phi_{T}):\phi,\phi_{t_{0}}=y\right\}}, \textrm{ as } \epsilon\downarrow 0.
\end{equation*}

A simple application of Jensen's inequality together with Varadhan's integral lemma (e.g., \cite{DemboZeitouni,FreidlinWentzell88,Varadhan}) shows that an asymptotically optimal $\bar{\mathbb{P}}$ should satisfy
\begin{equation*}
\lim_{\epsilon\rightarrow 0}\epsilon\ln\bar{\mathbb{E}}[e^{-\frac{1}{\epsilon}h(X^{\epsilon}_{T})}d\mathbb{P}/ d\bar{\mathbb{P}}]^{2}=-2G(t_{0},y), \textrm{ with } G(t,x)=\inf_{\phi\in\mathcal{AC}([t,T]), \phi_{t}=x}\left\{S_{tT}(\phi)+h(\phi_{T})\right\}\label{Eq:MinimumVariance3}
\end{equation*}
Turning to importance sampling, for $\bar{\mathbb{P}}$ that are absolutely continuous with respect to $\mathbb{P}$, Girsanov's formula implies %$\frac{d\bar{\mathbb{P}}}{d\mathbb{P}}=e^{-\frac{1}{2\epsilon}\int_{0}^{T}|v_{s}|^{2}
%ds+\frac{1}{\sqrt{\epsilon}}\int_{0}^{T}v_{s}dW_{s}},$
\begin{equation}
\frac{d\bar{\mathbb{P}}}{d\mathbb{P}}=e^{-\frac{1}{2\epsilon}\int_{0}^{T}|v_{s}|^{2}
ds+\frac{1}{\sqrt{\epsilon}}\int_{0}^{T}v_{s}dW_{s}}
\label{Eq:ChangeOfMeasure1}
\end{equation}
where $v_{t}$ is a progressively measurable
process (control) such that the right hand side  is a martingale (with respect to an appropriate filtration). Under $\bar{\mathbb{P}}$, $X^{\epsilon}$ satisfies %$X_{t_{0}}^{\epsilon }=y$ and
\begin{equation}
dX_{t}^{\epsilon }=\left[-\nabla V( X_{t}^{\epsilon })+\Gamma v_{t}\right]  dt+\sqrt{\epsilon }\Gamma d\bar{W}_{t},\quad \textrm{with}\quad\bar{W}_{t}=W_{t}-\frac{1}{\sqrt{\epsilon}}\int_{t_{0}}^{t}v_{\rho}d\rho\label{Eq:SimulationProcess}
\end{equation}

So, the problem is restricted to choosing the control $v_{t}$  optimally (i.e., such that the second moment is minimized)  and then using the estimator
based on iid samples generated from $\bar{\mathbb{P}}$ under (\ref{Eq:SimulationProcess}).
%So, the problem is restricted to choosing the control $v_{t}$  optimally (i.e., such that the second moment is minimized)  and then using the estimator
%(\ref{Eq:MonteCarloEstimate}), where $X_{k}$ are iid samples generated from $\bar{\mathbb{P}}$ under (\ref{Eq:SimulationProcess}).
%Let us define
%\begin{equation}
%G(t,x)= \inf_{\phi\in\mathcal{AC}([t,T]), \phi_{t}=x}\left\{S_{tT}(\phi)+h(\phi)\right\}.\label{Eq:MinimumVariance3}
%\end{equation}
Under appropriate conditions, the zero-variance (i.e. the best) change of measure is based on the control $v_{t}$ given by the formula $v_{t}=\bar{u}(t,X^{\epsilon}_{t})$ where
%\begin{equation}
$\bar{v}(t,x)=-\Gamma^{T} \nabla G^{\epsilon}(t,x) $ %\label{Eq:OptimalControlSmooth2}
%\end{equation}
where $G^{\epsilon}(t,x)$, with terminal condition $G^{\epsilon}(T,x)=h(x)$,  is the  solution to the PDE, of HJB type:
\begin{equation}
\partial_{t}G^{\epsilon}(t,x)-\nabla V(x)\cdot \nabla G^{\epsilon}(t,x)-\frac{1}{2}\left|\Gamma^{T} \nabla G^{\epsilon}(t,x)\right|^{2}+\frac{\epsilon}{2}\textrm{tr}\left[\Gamma\Gamma^{T} \nabla^{2}G^{\epsilon}(t,x)\right]=0.\label{Eq:HJBequationSmooth1}
\end{equation}

Since (\ref{Eq:HJBequationSmooth1}) is not tractable, it is standard approach to go to the viscosity limit $\epsilon\downarrow 0$. Then $G(t,x)=\lim_{\epsilon\downarrow 0}G^{\epsilon}(t,x)$  is the viscosity solution to the  HJB equation with Hamiltonian
\[
H(x,p)=\left<-\nabla V(x), p\right>-\frac{1}{2}\left\Vert \Gamma^{T} p\right\Vert^{2}
\]
 i.e., to the equation
     \begin{equation}
     \partial_{t}G(t,x)-\nabla V(x)\cdot DG(t,x)-\frac{1}{2}\left|\Gamma^{T} DG(t,x)\right|^{2}=0, \hspace{0.2cm} G(T,x)=h(x).\label{Eq:HJBequationSmooth}
     \end{equation}

Notice that by control arguments, e.g., see \cite{FleimingSoner2006}, we can also write
\[
G(t,x)=\lim_{\epsilon\downarrow 0}G^{\epsilon}(t,x)=\inf_{\phi\in\mathcal{AC}([t,T]), \phi_{t}=x}\left\{S_{tT}(\phi)+h(\phi_{T})\right\}.
\]

In fact, more is true. A smooth function $\bar{U}(t,x):[0,T]\times
\mathbb{R}^{d}\mapsto\mathbb{R}$ is called a \textit{subsolution} to the HJB
equation (\ref{Eq:HJBequationSmooth}) with $\epsilon=0$ if
\begin{equation}
    \partial_{t}\bar{U}(t,x)-\nabla V(x)\cdot \nabla\bar{U}(t,x)-\frac{1}{2}\left|\Gamma^{T} \nabla\bar{U}(t,x)\right|^{2}\geq 0, \hspace{0.2cm} \bar{U}(T,x)\leq h(x).\label{Eq:HJBequationSmoothSubsolution}
    \end{equation}

It turns out (Theorem 4.1 in \cite{DupuisSpiliopoulosWang}), that appropriate, \textit{smooth} subsolutions are enough.
If $\bar{U}(t,x)\in \mathcal{C}^{1,1}([t_{0},T]\times\mathbb{R}^{d})$ satisfies (\ref{Eq:HJBequationSmoothSubsolution}) and
the \emph{feedback} control to use in (\ref{Eq:SimulationProcess}) is $v_{t}=-\Gamma^{T} \nabla \bar{U}(t,X^{\epsilon}_{t})$, then
%\begin{equation}
%G(t_{0},y)+\bar{U}(t_{0},y)\leq \liminf_{\epsilon\rightarrow0}-\epsilon\ln \bar{\mathbb{E}}[e^{-\frac{1}{\epsilon}h(X^{\epsilon}_{T})}d\mathbb{P}/d\bar{\mathbb{P}}]^{2}\leq 2G(t_{0},y). \label{Eq:GoalRegime1Subsolution}%
%\end{equation}
\begin{equation}
G(t_{0},y)+\bar{U}(t_{0},y)\leq \liminf_{\epsilon\rightarrow0}-\epsilon\ln \bar{\mathbb{E}}\left[e^{-\frac{1}{\epsilon}h(X^{\epsilon}_{T})}\frac{d\mathbb{P}}{d\bar{\mathbb{P}}}\right]^{2}\leq 2G(t_{0},y). \label{Eq:GoalRegime1Subsolution}%
\end{equation}

Therefore,  asymptotic optimality is attained if $\bar{U}$ satisfies $\bar{U}(t_{0},y)=G(t_{0},y)=\lim_{\epsilon\downarrow 0}G^{\epsilon}(t_{0},y)$ since then lower and upper bound agree. The design and analysis of importance sampling schemes based on the systematic connection with subsolutions to the appropriate HJB and Isaacs  equations goes back to \cite{DupuisWang, DupuisWang2}. See also  \cite{BlanchetGlynn, BlanchetGlynnLeder2011, BlanchetGlynnLiu, BlanchetLiu2008} for the closely related concept of Lyapunov inequalities.% (essentially the exponentially of a subsolution).

The importance sampling simulation scheme in order to estimate $\theta^{\epsilon}(t_{0},y)\doteq\mathbb{E}_{t_{0},y}\left[  e^{-\frac{1}{\epsilon}h(X_{T}^{\epsilon}%
)}\right]  $ goes as follows. Let $X^{\epsilon,v}$ be
the solution to the SDE
\begin{align}
dX_{t}^{\epsilon,v}&=\left(-\nabla V( X_{t}^{\epsilon,v})+\Gamma v_{t} \right) dt+\sqrt{\epsilon }\Gamma dW_{t},\hspace{%
0.2cm}X_{t_{0}}^{\epsilon,u}=y.\label{Eq:LDPandA2}
\end{align}

\begin{enumerate}
\item{Consider $v_{t}=\bar{u}(t,X_{t}^{\epsilon,v})=-\Gamma^{T}\nabla_{x}\bar{U}(t,X_{t}^{\epsilon,v})$  with $\bar{U}$ an appropriate subsolution, i.e., it satisfies (\ref{Eq:HJBequationSmoothSubsolution})}
\item {Consider the estimator
\begin{equation}
\hat{\theta}^{\epsilon}(y)\doteq\frac{1}{N}\sum_{j=1}^{N}\left[  e^{-\frac
{1}{\epsilon}h(X_{T}^{\epsilon,v}(j))}Z_{j}^{v}\right]
\label{Def:OptimalEstimator1}%
\end{equation}
where
\[
Z_{j}^{v}\doteq e^{-\frac{1}{2\epsilon}\int_{0}^{T}\left\| \bar{u}\left(
t,X_{t}^{\epsilon,v}(j)\right)  \right\| ^{2}dt-\frac{1}{\sqrt{\epsilon}%
}\int_{0}^{T}\left<\bar{u}\left(  t,X_{t}^{\epsilon,v}(j)\right),  dW_{t}%
(j)\right>}\label{Def:OptimalChangeOfMeasure1}%
\]
and $(W(j),X^{\epsilon,v}(j))$ is an independent sample generated from
(\ref{Eq:LDPandA2}) with control $v_{t}=\bar{u}\left(  t,X_{t}^{\epsilon
,v}(j)\right)  $.}
\end{enumerate}

We conclude this section, with the remark that a choice of the control $v_{t}$ based on a subsolution as defined by (\ref{Eq:HJBequationSmoothSubsolution})  only guarantees logarithmic asymptotic optimality and does not say something about the \textit{important} effect of pre-factors. As we will see in Section \ref{S:EffectRestPoints}, this can imply degradation in the performance of the algorithm in problems with metastability.  When dealing with metastability issues, things may be even more problematic if one is using the exact solution to the association HJB equation, $G(t,x)$. While this may be not be a problem for problems that do not involve rest points (i.e. does not involve stable or unstable equilibrium points) in the domain of interest, it does become problematic when dealing with metastability issues.

\begin{remark}
Obtaining  accurately the solution $G(t,x)$ to the HJB equation (\ref{Eq:HJBequationSmooth}), analytical or numerical,  is challenging in high dimensions. However, even if this were possible, the solution by itself is not always  suitable for importance sampling when one is interested in computing escape or transition probabilities.
The issue is that in these cases, the solution is a viscosity solution with a discontinuous derivative at the rest point (stable or unstable equilibrium points) and with negative definite generalized second derivative there. Physically, the exact solution to the HJB equation attempts at each point in time and space to force the simulated trajectories to follow a most likely large deviations optimal path. However, by standard control arguments, see \cite{FleimingSoner2006}, the discontinuity of the spatial derivative at the rest point, implies that multiple most likely optimal paths exist. As a consequence, the noise
can cause trajectories to return to a neighborhood of the origin, thereby
producing large likelihood ratios.  In Section \ref{SS:IS_MostLikelyPaths}, we will see that this is a serious issue, leading to poor performance, even in dimension one where one can solve the HJB equation analytically. Importance sampling, when dealing with state dependent metastable dynamical systems, needs to be addressed from a global point of view and not local.
\end{remark}

\section{The effect of rest points  on importance sampling}\label{S:EffectRestPoints}

%In \cite{DupuisSpiliopoulosZhou2013,DupuisSpiliopoulos2014a,DupuisSpiliopoulos2014b,Spiliopoulos2013b} efficient Monte Carlo schemes for estimating the probability
%cumulative distribution of exit times in the presence of rest points are constructed for processes of the form (\ref{Eq:LDPnoHomogenization}).

As it is shown, mathematically and numerically, in \cite{DupuisSpiliopoulosZhou2013,DupuisSpiliopoulos2014a,Spiliopoulos2013b}, in dynamical systems that exhibit metastable behavior standard simulation methods do not readily apply. Asymptotic optimality is \textit{necessary but not sufficient} for good performance
due to the non-trivial effect of the pre-factors.  The pre-factor computations in \cite{DupuisSpiliopoulosZhou2013,Spiliopoulos2013b} prove that there is non-trivial interaction of parameters such as
the strength of the noise $\epsilon$ and the terminal time $T$. We remark here that this is in contrast to
escape probabilities for other well studied problems, such as stochastic networks, e.g., \cite{BlanchetGlynn, BlanchetGlynnLiu, DupuisLederWang2007,DupuisLederWang2009,DupuisSezerWang2007,DupuisWang2},
because there the proximity of the rest point has little impact on either the asymptotic rate of decay or the pre-exponential term.

These interactive effects vanish in the logarithmic limit as the noise goes to zero, but they have
a significant effect on the performance of the algorithms. The following question immediately presents itself:
\begin{itemize}
\item{Is it sufficient to have schemes that are only asymptotically logarithmical optimal, in the sense that the second moment of the estimator satisfies (\ref{Eq:GoalRegime1Subsolution})? What about pre-factors? Are they truly negligible in practice in the rare event regime?}
\item{Can we construct a subsolution $\bar{U}(t,x)$ that not only satisfies (\ref{Eq:HJBequationSmoothSubsolution}) but it also takes care of the prefactor effects?}
\end{itemize}

\subsection{Effects in the prelimit}\label{SS:QuasipotentialSubsolution}
Let us demonstrate the effect of prefactors on the behavior of estimators in the following classical simple setting. Let us assume that the diffusion coefficient $\Gamma=I$, and that $x=O$ is the global minimum for $V(x)$. In particular, let us assume that $DV(O)=0$ and that $DV(x)\neq 0$ for every $x\neq O$. Define
\[
\mathcal{D}=\left\{ x\in\mathbb{R}^{d}: 0\leq V(x)<L \right\}
\]
and let $A_{c}=\left\{ x\in\mathbb{R}^{d}: V(x)=c\right\}$. Then for an initial point $y$ such that $0\leq V(y)<L$, let us assume that we want to estimate
\[
\theta^{\epsilon}(t,y)=\mathbb{P}_{t,y}\left\{ X^{\epsilon} \text{ hits } A_{L} \text{ before time }T\right\}.
\]

A classical quantity if interest in metastability theory is the quasipotential, see \cite{FreidlinWentzell88}. The quasipotential with respect to the equilibrium point $O$ is defined as follows
\[
W(O,x)=\left\{ S_{0T}(\phi): \phi\in \mathcal{C}([0,T]), \phi(O)=0, \phi(T)=x, T\in(0,\infty)\right\}
\]

Under the assumptions of this section, the  quasi-potential is computable in closed form \cite{FreidlinWentzell88}:
%\begin{equation*}
$W (O,x) = 2V(x) \textrm{ for } x\in\{y\in \mathcal{D}\cap\partial \mathcal{D}: V(y)\leq \inf_{z\in\partial \mathcal{D}}V(z)\}.$
%\end{equation*}
%Additionally, the unique extremal of the functional $S(\phi)$ on the set of functions
%$\phi_{s}, -\infty\leq s\leq T$ leading from $O$ to $x$ is given by the equation
%$\dot{\phi}_{s}=DV(\phi_{s}), s\in(-\infty,T), \phi_{T}=x.$

Now, if we define $\tau^{\epsilon}=\inf\left\{t>0: X_{t}^{\epsilon}\notin \mathcal{D}   \right\}$, then, as it is shown in \cite{FreidlinWentzell88} we have that
%\begin{equation*}
 $\lim_{\epsilon\downarrow 0}\epsilon\ln\mathbb{E}\tau^{\epsilon}=\inf_{z\in\partial \mathcal{D}}W(O,z).$
%\end{equation*}
Thus, the quasi-potential allows to approximate exit times in the logarithmic large deviations regime,  \cite{FreidlinWentzell88}.  Many quantities in the theory of metastability are defined via the quasi-potential.
The quasi-potential characterizes the leading asymptotics of exit times and exit probabilities, approximates transition rates for reversible and irreversible systems and allows to qualitatively describe transitions between stable attractors if the system has many of them; see also \cite{Cameron2012,DayDarden1985,FlemingJames1992,FreidlinWentzell88,MaierStein1993,MaierStein1997} for more details. These conclusions hold for both gradient and non-gradient cases, but in the gradient case the quasi-potential is computable in closed form.

Turning now to importance sampling,  it is easy to verify that the quasi-potential is a stationary subsolution to the associated HJB equation (\ref{Eq:HJBequationSmoothSubsolution}) with $\epsilon=0$, by adding an appropriate constant $C$ in order to justify the necessary boundary and terminal conditions. In particular, $\bar{U}_{QP}(x)=2L-W(O,x)$ defines a subsolution for (\ref{Eq:HJBequationSmoothSubsolution}). It turns out, see \cite{DupuisSpiliopoulosZhou2013}, that the quasipotential yields a  reasonable change of measure if rest points are not part of the domain of interest. However, this is no longer true if rest points are included in the domain of interest.

Let us denote $Q^{\epsilon}(0,y;\bar{u})=\bar{\mathbb{E}}[e^{-\frac{1}{\epsilon}h(X^{\epsilon}_{T})}d\mathbb{P}/ d\bar{\mathbb{P}}]^{2}$ to be the second moment of the estimator constructed using the control $\bar{u}$. Based now on the arguments of \cite{DupuisSpiliopoulosZhou2013} one can prove the following representation for the second moment of the estimator estimator based on the change of measure induced by the control $\bar{u}(t,x)=-\nabla \bar{U}_{QP}(x)$
\begin{equation}
-\epsilon\log Q^{\epsilon}(0,y;\bar{u})=\inf_{v\in\mathcal{A}}%
\mathbb{E}\left[  \frac{1}{2}\int_{0}^{\hat{\tau}^{\epsilon}}\left\Vert
v(s)\right\Vert ^{2}ds-\int_{0}^{\hat{\tau}^{\epsilon}}\left\Vert \bar
{u}(\hat{X}^{\epsilon}_{s})\right\Vert ^{2}ds+\infty1_{\left\{  \hat{\tau
}^{\epsilon}>T\right\}  }\right]  . \label{Eq:game_rep}%
\end{equation}
 where $\hat{X}^{\epsilon}_{s}$ is the unique solution to the SDE
 \[
d\hat{X}^{\epsilon}_{s}=-DV(\hat{X}^{\epsilon}_{s})ds+\left[
\sqrt{\epsilon}dW_{s}-[\bar{u}(\hat{X}^{\epsilon}_{s})-v(s)]ds\right]
\label{Eq:control_diffusion}%
\]
with initial condition $\hat{X}^{\epsilon}_{0}=y$ and $\hat{\tau
}^{\epsilon}$ is the first time that $\hat{X}^{\epsilon}$ exits from $\mathcal{D}$.

It is important to note that (\ref{Eq:game_rep}) provides a non-asymptotic representation for the second moment of the estimator. By the arguments of \cite{DupuisSpiliopoulosZhou2013}, we can choose a particular admissible control $v(s)$ in (\ref{Eq:game_rep}) so that the following takes place. Let $T$ be large and let $0<K<T$ so that the time interval $[0,T]$ is split into $[0,T-K)$ and $[T-K,T]$. Set $v(s)=0$ for $s\in[0,T-K)$. The resulting dynamics for $\hat{X}^{\epsilon}$ is stable for $s\in[T-K,T]$ and with high probability the process will stay around the point $y$ for $s\in[0,T-K)$. In the time interval $[T-K,T]$, we set $v(s)$ so that escape happens prior to $T$. Then, it can be shown that there are positive constants $C_{1},C_{2}<\infty$, so that
\[
Q^{\epsilon}(0,y;\bar{u})\geq e^{-\frac{1}{\epsilon}C_{1}+C_{2}(T-K)}.
\]

This bound indicates that if $T$ is large, one may need to go to considerably small values of $\epsilon$ in order to achieve the theoretical optimal asymptotic performance. We also remark that if $T$ is large (see Chapter 4 of \cite{FreidlinWentzell88}), $G(0,y)$ and $\bar{U}(y)$ get closer in value. Thus, by (\ref{Eq:GoalRegime1Subsolution}) and for large enough $T$, the particular importance sampling scheme is asymptotically optimal.

Hence, we have just seen an example where an importance sampling estimator is almost asymptotically optimal, but it does not perform that well pre-asymptotically due to the effect of the possibly long time horizon $T$ and its interplay with $\epsilon$.

\subsection{The problems arising when following large deviations asymptotically most likely paths and a remedy to the problem}\label{SS:IS_MostLikelyPaths}
The connection of change of measures with HJB equations via large deviations is well situated for a systematic treatment of dynamic importance sampling schemes  for state dependent processes like diffusions (\ref{Eq:LDPnoHomogenization}). For small noise diffusions the theoretical framework of subsolutions to HJB equations and their use for Monte Carlo methods can be found in \cite{DupuisSpiliopoulosWang}. It was a common belief for sometime that if the underlying stochastic process has a large deviations principle and if the change of measure is consistent with
 the large deviations asymptotically most likely path leading to the rare event (an open-loop control), then the resulting importance sampling scheme would be optimal. %However, this
%heuristic has been theoretically justified only for very special systems.
However, such heuristics have been shown to be unreliable in general and simple examples have been
constructed showing the failure of the corresponding importance schemes even in very simple settings \cite{GlassermanWang1997, GlassermanKou1995}. This is due to the presence
of ``rogue-trajectories", i.e., unlikely trajectories, that are likely enough to increase likelihood ratios to the point that
 the performance is comparable to standard Monte Carlo. This is especially true for metastability problems (i.e., when transitions between fixed points occur at suitable (large) timescales)
 where multiple nearly optimal paths may exist.

 Use of dynamic changes of measure, i.e. based on feedback controls (time and location dependent) becomes important, see \cite{DupuisSpiliopoulosZhou2013,DupuisSpiliopoulos2014a}. However, even changes of measures that are based on feedback controls, that are consistent with large deviations and lead to asymptotically optimal change of measures can also be problematic in practice. We demonstrate this below in Table \ref{T:Table1}. Namely, as it turns out, in the presence of rest points and metastability, the prefactors may affect negatively the behavior of estimators even if one is using asymptotically optimal changes of measure in the spirit of (\ref{Eq:GoalRegime1Subsolution}). Hence, it becomes important to use dynamic change of measures that are based on subsolutions but lead to good performance even pre-asymptotically.

To that end, \textit{novel} explicit simulation schemes are then constructed in \cite{DupuisSpiliopoulosZhou2013,DupuisSpiliopoulos2014a} that perform provably-well
\textit{both} asymptotically and non-asymptotically, even when the  simulation time is long.
These constructions are based on large deviations asymptotics \cite{BovierEckhoffGayrandKlein2004,BovierGayrandKlein2005,FreidlinWentzell88}, stochastic control arguments and asymptotic expansions \cite{FlemingJames1992, FleimingSoner2006}  and detailed asymptotic analysis of the subsolution to the
associated HJB in the neighborhood of the rest point where the potential can be thought of as being approximately quadratic. Essentially, due to the fact that near the rest point, the potential can be thought of as being approximately quadratic, one can hope to solve or to approximate the solution to the associated variational problem there. Then one needs to patch this solution together with the quasipotential based subsolution (which is a good subsolution away from the rest point) in the right way. Then, the combined subsolution, see $\bar{U}^{\delta}(t,x)$ in (\ref{Eq:ChangeOfMeasure}), turns out to be a good approximation to the zero variance change of measure. Such schemes lead to importance sampling algorithms with provably-good performance for all small $\epsilon>0$ and without suffering from bad prefactor effects.

In order to illustrate the point, let us briefly demonstrate such a construction in the case of dimension one, see \cite{DupuisSpiliopoulosZhou2013}. So, let us assume that $V(x)=\frac{\lambda}{2}x^{2}$ with $\lambda>0$ and let us assume that we study the problem of crossing a level set, say $L$, of the potential function $V(x)$. Here, we can compute $G(t,x)$ in closed form and we get
\begin{align}
G(t,x) &  =\inf_{\phi_{t}=x,V(\phi_{T})=L}\left\{  \frac{1}{2}\int_{t}%
^{T}\left\Vert \dot{\phi}_{s}+\lambda\phi_{s}\right\Vert ^{2}ds\right\}
  =\inf_{\hat{x}\in V^{-1}(L)}\lambda\frac{\left(  \hat{x}-xe^{\lambda
(t-T)}\right)  ^{2}}{1-e^{2\lambda(t-T)}}.\label{Eq:VariationalProblem}%
\end{align}

Notice, that $G(t,x)$ is also a viscosity solution to the $\epsilon=0$ HJB equation  (\ref{Eq:HJBequationSmooth}) when supplemented with the appropriate boundary conditions. Hence, based on (\ref{Eq:GoalRegime1Subsolution}) a change of measure based on $G(t,x)$, i.e., using the control $u(t,x)=-\partial_{x}G(t,x)$, is expected to yield an asymptotically efficient estimator. While this is true, we will see below that this is not sufficient to yield good performance. The fact that the function $G(t,x)$ is not continuously differentiable in the domain of interest, implies that multiple optimal paths exist, which is an intuitive reason for the degradation in performance that will be demonstrated below.

However, by appropriately mollifying  $G(t,x)$ and combining it with the quasipotential subsolution (as constructed in Section \ref{SS:QuasipotentialSubsolution}), one can construct a global subsolution which performs provably well even pre-asymptotically. The point  is that $G(t,x)$ provides a good change of measure
while near the rest point, whereas the quasipotential induced subsolution $\bar{U}_{QP}(x)=2L-W(O,x)$ provides a
good change of measure away from the rest point. There are a few more issues
to deal with though. The first one is that $G(t,x)$ is discontinuous near
$t=T$. The second one is that we need to put them together in a smooth way
that will define a global subsolution.

Since $G(t,x)$ is discontinuous at $t=T$,  we introduce two
mollification parameters $t^{\ast}$ and $M$ that will be appropriately chosen
as functions of $\epsilon$. Motivated by the fact that $G(t,x)$ is a good subsolution near the equilibrium point, we fix another parameter
$\hat{L}\in(0,L]$. In the one-dimensional case, it is easy to solve the
equation $V(x^{*})=\hat{L}$ and in particular we get that $x^{*}=\pm \hat{x}$ where $\hat{x}%
=\sqrt{\frac{2\hat{L}}{\lambda}}$. As a matter of fact, instead of using $G(t,x)$ directly,
we set
\begin{align*}
F^{M}(t,x;\hat{x})  &  = \lambda\frac{\left(  \hat{x}-xe^{\lambda
(t-T)}\right)  ^{2}}{\frac{1}{M}+1-e^{2\lambda(t-T)}}%
\end{align*}

In order now to pass smoothly between the $\bar{U}_{QP}(x)$ and $F^{M}(t,x;\hat{x})$ or $F^{M}(t,x;-\hat{x})$ without violating the subsolution property, we use the
exponential mollification, see \cite{DupuisWang2}
\[
U^{\delta}(t,x)=-\delta\log\left(  e^{-\frac{1}{\delta}\bar{U}_{QP}(x)}+e^{-\frac
{1}{\delta}\left[  F^{M}(t,x;\hat{x})+\bar{U}_{QP}(\hat{x})\right]  }+e^{-\frac
{1}{\delta}\left[  F^{M}(t,x;-\hat{x})+\bar{U}_{QP}(-\hat{x})\right]  }\right)
\]
It is easy to see that as $\delta\downarrow 0$
\[
\lim_{\delta\downarrow 0}U^{\delta}(t,x)=\min\{\bar{U}_{QP}(x), F^{M}(t,x;\hat{x}), F^{M}(t,x;-\hat{x})\}
\]

Clearly, if we choose $\hat{L}=L$, then we get $\bar{U}_{QP}(\hat{x})=0$. Based on these constructions, a provably efficient importance sampling scheme is constructed in
\cite{DupuisSpiliopoulosZhou2013}, based on the
subsolution
\begin{align}
\bar{U}^{\delta}(t,x)=\left\{
\begin{array}
[c]{cc}%
\bar{U}_{QP}(x), & t>T-t^{\ast}\\
U^{\delta}(t,x), & t\leq T-t^{\ast}%
\end{array}
\right.\label{Eq:ChangeOfMeasure}
\end{align}

It turns out that $\bar{U}^{\delta}(t,x)$ is a global smooth subsolution which
has provably good performance both pre-asymptotically and asymptotically.  The role of the exponential mollification is to allow a smooth transition between the region that is near the equilibrium point and the region that is far away from it. %The
The
precise optimality bound and its proof guide the choice of the parameters
$\delta,t^{\ast},M$ and $\hat{L}$.
For the convenience of the reader, we present in Table \ref{TableParameters} the suggested values for $(\delta,\hat{L},M,t^{*})$, given the value of
the strength of the noise $\epsilon>0$.

 \begin{table}[th]
\begin{center}
{\small
\begin{tabular}
[c]{|c|c|c|c|c|}\hline
parameter & $\delta$ & $\hat{L}\in(0,L]$ & $M$ & $t^{*}$ \\\hline
values & $2\epsilon$ & $O(1)$ or $\varepsilon^{2m}$ with $m<\kappa$ & $\max\{\frac{\hat{L}}{\varepsilon^{2\kappa}}, 4\}$ with $\kappa\in(0,1/2)$ & $-\frac{2}{\lambda}\log\frac{1}{M}$ \\\hline
\end{tabular}}
\end{center}
\caption{Parameter values for the algorithm based on a given value of $\epsilon>0$.}\label{TableParameters}
\end{table}

We refer the interested reader to \cite{DupuisSpiliopoulosZhou2013,DupuisSpiliopoulos2014a} for further details on the theoretical performance of the algorithm and on the choice of parameters.

In order to illustrate in a simple setting the effect of prefactors in the presence of metastable effects, we record in Table \ref{T:Table1} Monte Carlo estimates based on $K=10^{7}$ trajectories for the exit time distribution $\mathbb{P}_{y}\left[\tau^{\epsilon}_{\mathcal{D}\cup \partial \mathcal{D}}\leq T\right]$  from the basin of attraction of the left attractor of
the potential of Figure 1 for the process $X^{\epsilon}$ given by (\ref{Eq:LDPnoHomogenization}) with $\Gamma=I$. We used the importance sampling (IS) methods of \cite{DupuisSpiliopoulosZhou2013}, i.e., the change of measure based on the subsolution (\ref{Eq:ChangeOfMeasure}) and record estimates for different pairs $(\epsilon,T)$.
%The table records estimated probabilities for
%different pairs $(\epsilon,T)$, using the optimal (both pre-asymptotically and asymptotically) change of measure constructed in \cite{DupuisSpiliopoulosZhou2013}.
In the figures next to Table \ref{T:Table1}, we  compare the relative errors per sample of (a): the algorithm, which is optimal for all $\epsilon>0$, i.e the one based on the subsolution $\bar{U}^{\delta}(t,x)$, with (b): the IS algorithm that is consistent with the
large deviations asymptotically most likely path leading to the rare event, i.e the one based on the actual solution $G(t,x)$ of the associated HJB equation. Notice  however that the IS algorithm based on $G(t,x)$ is only asymptotically optimal in the large deviations logarithmic sense as $\epsilon\downarrow 0$ (i.e., it satisfies (\ref{Eq:GoalRegime1Subsolution})).

\begin{table}[th]
{\scriptsize
\begin{minipage}{.55\textwidth}
\centering
\begin{tabular}
[c]{|c|c|c|c|c|c|}\hline
$\epsilon\hspace{0.1cm} | \hspace{0.1cm} T$ &  $2.5$  & $7$ &$10$  & $18$ & $23$\\\hline
$0.20$  & $2e-02$  & $8.3e-02$ & $1.2e-01$ & $2.1e-01$ & $2.7e-01$\\\hline
$0.16$  & $7e-03$  & $2.7e-02$ & $4.0e-02$ & $7.4e-02$ & $9.5e-02$\\\hline
$0.13$  & $2e-03$  & $6.9e-03$ & $1.1e-02$ & $2.0e-02$ & $2.6e-02$\\\hline
$0.11$  & $4e-04$  & $1.8e-03$ & $2.8e-03$ & $5.4e-03$ & $7.0e-03$\\\hline
$0.09$  & $5e-05$  & $2.6e-04$ & $4.1e-04$ & $7.8e-04$ & $1.0e-03$\\\hline
$0.07$  & $2e-06$  & $1.2e-05$ & $1.9e-05$ & $3.7e-05$ & $4.8e-05$\\\hline
$0.05$  & $7e-09$  & $4.4e-08$ & $7.0e-08$ & $1.4e-07$ & $1.8e-07$\\\hline
\end{tabular}
% \caption{Exit time probability estimates $\mathbb{P}_{y}\left[\tau^{\epsilon}_{\mathcal{D}\cup \partial \mathcal{D}}\leq T\right]$  for different pairs $(\epsilon,T)$, using the optimal change of measure constructed in \cite{DupuisSpiliopoulosZhou2013}.
% \label{T:Table1}
\end{minipage}
\begin{minipage}{.45\textwidth}
\centering
\includegraphics[width=0.60\linewidth]{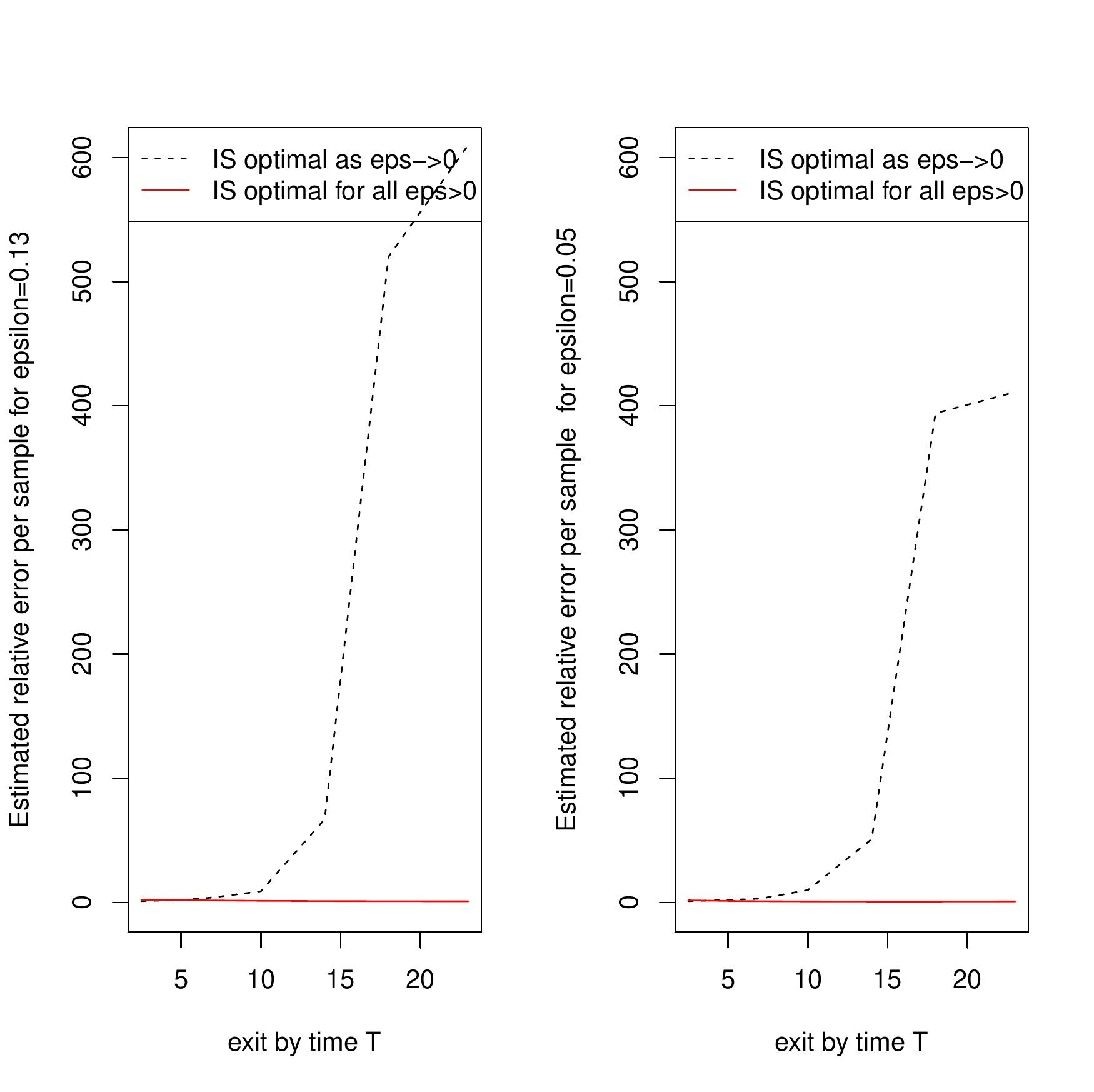}
%\caption{Comparison of relative errors per sample for exit probabilities by time T.}\label{F:Figure1}
\end{minipage}
}
\caption{Left: Exit time distribution $\mathbb{P}_{y}\left[\tau^{\epsilon}_{\mathcal{D}\cup \partial \mathcal{D}}\leq T\right]$  for different pairs $(\epsilon,T)$,
using the optimal change of measure constructed in \cite{DupuisSpiliopoulosZhou2013}. Events range from very rare to not so rare.
Right: Comparison of relative errors per sample for two different changes of measure and for two values of $\epsilon$. Small relative error is better. }\label{T:Table1}
%\caption{Left: exit time probability estimates $\mathbb{P}_{y}\left[\tau^{\epsilon}_{\mathcal{D}\cup \partial \mathcal{D}}\leq T\right]$  for different pairs $(\epsilon,T)$,
%using the optimal change of measure constructed in \cite{DupuisSpiliopoulosZhou2013}.
%Right: Comparison of relative errors per sample for two differnt changes of measure, one that is only asymptotically optimal and one that is also
%pre-asymptotically optimal based on \cite{DupuisSpiliopoulosZhou2013}. The smaller the relative error is, the better it is. }\label{T:Table1}
\end{table}

Using relative error per sample as comparison criterium, we compare the two algorithms for two values of $\epsilon$, one for which the events are not so rare ($\epsilon=0.13$) and one for which the events are very rare ($\epsilon=0.05$). Exact values are in the table, and we remark for completeness that intermediate behavior is qualitatively the same. Both algorithms perform well when $T$ is small, but the algorithm that is based on the solution of the associated HJB equation, which is only logarithmic asymptotically optimal,
starts deteriorating considerably as $T$ gets large. The latter  is an effect of the pre-factors becoming important. On the other hand, the change of measure
constructed in \cite{DupuisSpiliopoulosZhou2013} that takes into account the pre-factor information and is pre-asymptotically optimal, yields optimal performance independently of the values $\epsilon$ and $T$ with relative errors around one, meaning that the values recorded at the table are reliable. It is important to note that due to large deviations, exit happens in long time scales, which implies that reliable estimates, especially when $T$ is large, are essential.

\section{Importance sampling for rough energy landscapes}\label{S:IS_MultiscaleLangevin}
In Section \ref{S:EffectRestPoints}, we reviewed some of the practical issues that come up when  one is trying to apply importance sampling techniques to metastable dynamics. While in Section \ref{S:EffectRestPoints} we ignored the effect of multiple scales, the goal of this section is to address the role of multiple scales in the design of asymptotically optimal importance sampling schemes.

 A particular model of interest in chemical physics is the first order Langevin equation (\ref{Eq:LangevinEquation2}).  Let us consider
\begin{equation}
dX_{t}^{\epsilon,\delta }=\left[ -\frac{\epsilon}{\delta }\nabla Q\left(
X_{t}^{\epsilon,\delta }/\delta \right) -\nabla V\left( X_{t}^{\epsilon,\delta }\right) %
\right] dt+\sqrt{\epsilon }\sqrt{2D}dW_{t},\hspace{0.2cm}X_{0}^{\epsilon,\delta
}=y, \hspace{0.2cm}, 0<\epsilon,\delta\ll 1,  \label{Eq:LangevinEquation2}
\end{equation}
where  the two-scale potential is composed
by a large-scale part, $V(x)$, and a fluctuating part, $\epsilon Q(x/\delta)$. If $Q$ is periodic then we have a periodic environment, whereas if $Q$ is random then we have a random environment. Models like (\ref{Eq:LangevinEquation2}) can be used to model rough energy landscapes
\cite{BanushkinaMeuwly2007,Janke,Zwanzig,DupuisSpiliopoulosWang2}.  As it has been suggested
(e.g., \cite{LifsonJackson, Zwanzig}), the associated energy landscapes of certain
biomolecules can be rugged (i.e., consist of many local ``small" minima within local deep minima separated
by barriers of varying heights). When one is interested in rare events, large deviations and Monte Carlo methods are relevant.

\begin{figure}[ht]\label{F:Figure1}
\begin{center}
\includegraphics[scale=0.3, width=6.0 cm, height=2.0 cm]{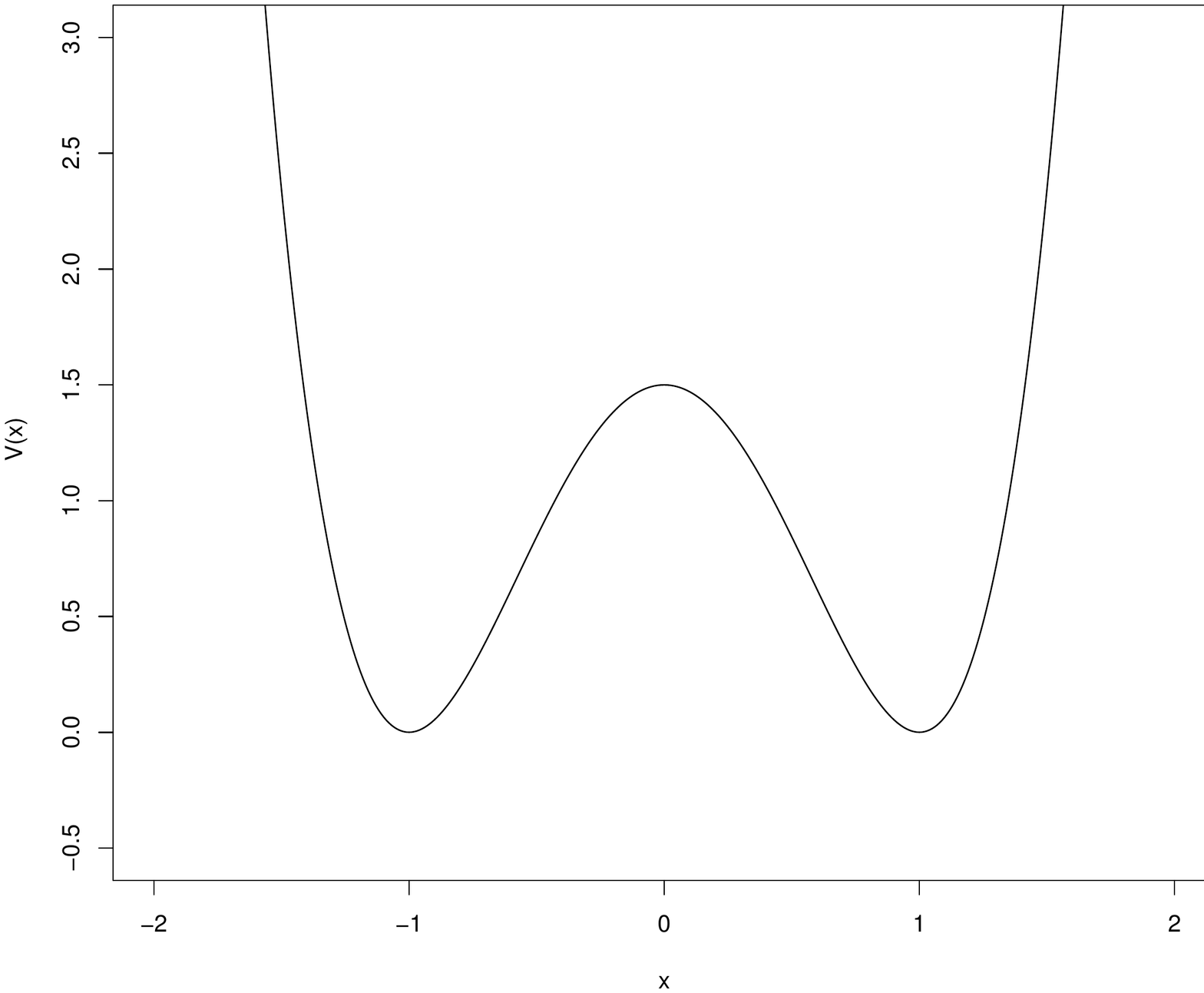}\hspace{2cm}
\includegraphics[scale=0.3, width=6.0 cm, height=2.0 cm]{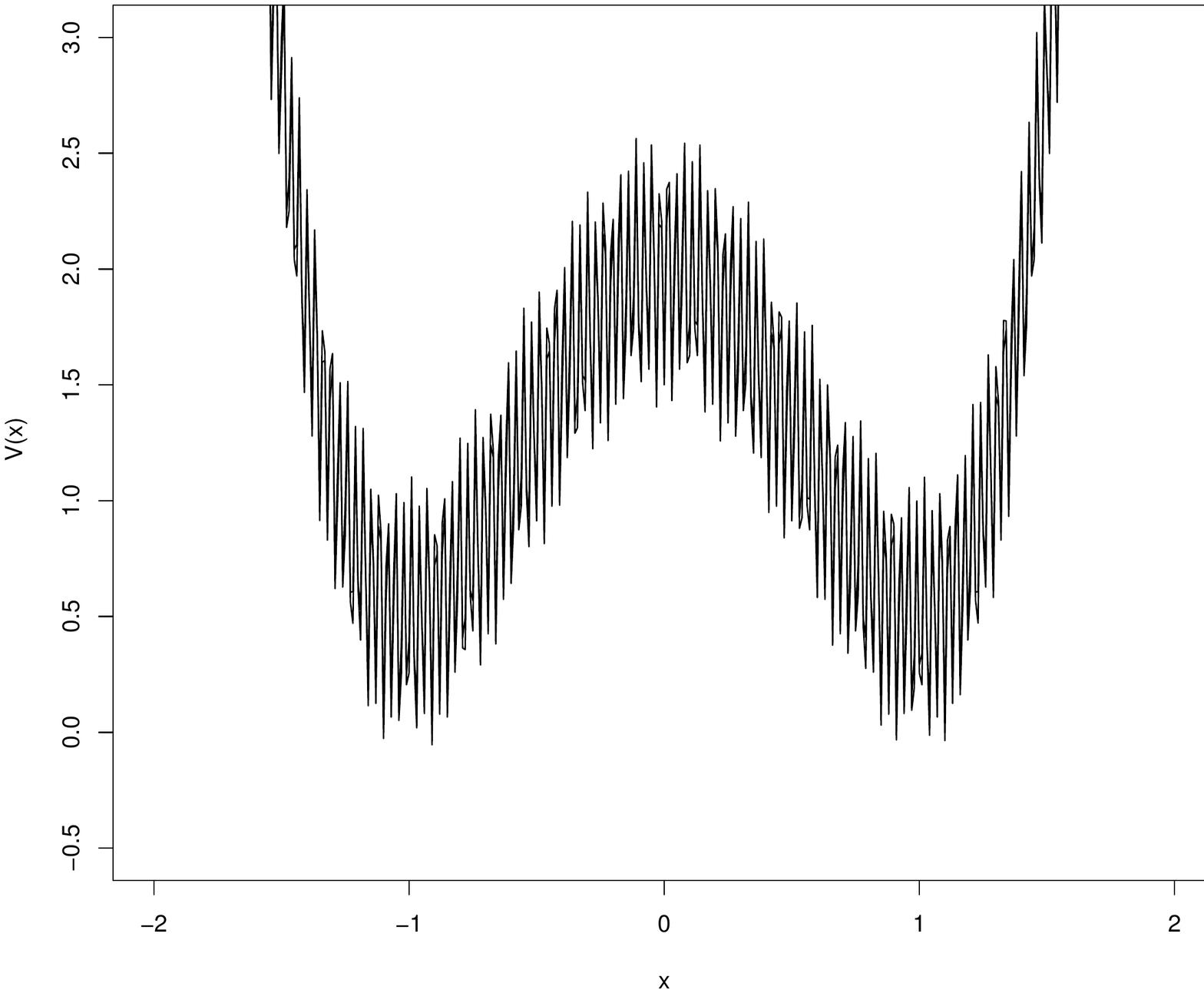}
\caption{A smooth and a rough potential function (energy landscape) with two wells.}
\end{center}
\end{figure}

If $Q(y)$ is {\it periodic},  large deviations for multiscale diffusions in periodic environments are obtained in \cite{DupuisSpiliopoulos,FS,Spiliopoulos2013a} for all possible interactions between $\epsilon$ and $\delta$, setting the ground for the mathematical formulation of the related importance sampling theory, \cite{DupuisSpiliopoulosWang,DupuisSpiliopoulosWang2,Spiliopoulos2013a}. The {\em novel feature} is that the optimal change of measure for importance sampling is not based only on the gradient of the homogenized HJB equation (as in Subsection \ref{S:ReviewIS_LDP}). The effect of fluctuations, which is quantified via the solution to the ``cell problem" in  homogenization \cite{BLP,PS},
is equally important. The cell problem is the solution to a Poisson type PDE. It is used to define the so called ``corrector", which characterizes the first order correction in the approximation of the multiscale HJB  by its homogenized limit.  Therefore, when compared to the case without multiple scales, one needs more detailed information in order to guarantee, \textit{at least}, asymptotic optimality.

For example, consider model (\ref{Eq:LangevinEquation2}) in the case $\frac{\epsilon}{\delta}\uparrow\infty$. Define the Gibbs measure
\[
\mu(dy)=\frac{1}{L} e^{-\frac{Q(y)}{D}}dy,\quad L=\int_{\mathbb{T}^{d}} e^{-\frac{Q(y)}{D}}dy.
\]

Then denote by $\chi(y)$  the smooth solution to the ``cell problem''
\begin{equation}
-\nabla Q(y)\cdot \nabla \chi(y)+ D\textrm{tr}\left[\nabla^{2}\chi(y)\right] = \nabla Q(y),\qquad \int \chi(y)\mu(dy) = 0.
\end{equation}

The following large deviations result holds which is a special case of the results of \cite{DupuisSpiliopoulos}. In particular, \cite{DupuisSpiliopoulos} covers the case of general state dependent drift (not necessarily of gradient form) and state dependent diffusion coefficient.
\begin{theorem}[Theorem 5.3 of \cite{DupuisSpiliopoulos} for the case of (\ref{Eq:LangevinEquation2})]
\label{T:MainTheorem3} Assume that the functions $\nabla Q\left(
y \right)$ and $\nabla V\left(x\right)$ are continuous and globally bounded,
as are their partial derivatives up to order $1 $ in $y$ and order $2$ in $x$ respectively. Let $\{X^{\epsilon,\delta},\epsilon,\delta>0\}$ be the unique
strong solution to (\ref{Eq:LangevinEquation2}). Let
\begin{align*}
r(x)  &  =-\int_{\mathbb{T}^{d}}\left(  I+\frac{\partial\chi(y)}{\partial
y}\right)  \mu(dy) \nabla V(x),\\
q  &  =2D \int_{\mathbb{T}^{d}}\left(  I+\frac{\partial\chi(y)}{\partial
y}\right)  \left(  I+\frac{\partial\chi(y)
}{\partial y}\right)^{T}\mu(dy),
\end{align*}
where $I$ denotes the identity matrix. If $\epsilon/\delta\rightarrow \infty$, then $\{X^{\epsilon,\delta},\epsilon,\delta>0\}$ converges in probability as $\epsilon,\delta\rightarrow 0$ to the solution of the ODE
\[
d\bar{X}_{t}=r(\bar{X}_{t})dt
\]
and satisfies a large deviations principle with rate function
\[
S_{tT}(\phi)=%
\begin{cases}
\displaystyle{\frac{1}{2}\int_{t}^{T}\left(  \dot{\phi}_{s}-r(\phi_{s})\right) q^{-1}\left(  \dot{\phi}_{s}-r(\phi_{s})\right)^{T} ds } & \text{if }\phi\in\mathcal{AC}%
([t,T]),\phi_{t}=x\\
+\infty & \text{otherwise}.
\end{cases}
\]
\end{theorem}

In addition, it turns out that an asymptotically efficient change of simulation measure can be constructed analogously to Section \ref{S:EffectRestPoints}, but  based on the  feedback control (see Theorem 4.1 in \cite{DupuisSpiliopoulosWang})
\begin{equation}
v_{t}=\bar{u}(t,X^{\epsilon}_{t},X^{\epsilon}_{t}/\delta),\quad \textrm{with}\quad \bar{u}(t,x,y)=-\sqrt{2D}\left( I+ \partial \chi(y)/\partial y\right) ^{T}\nabla_{x} \bar{U}(t,x).\label{Eq:OptimalControlRough}
\end{equation}

$\bar{U}(t,x)$ satisfies the inequalities in (\ref{Eq:HJBequationSmoothSubsolution}) with the homogenized (averaged) coefficients $r(x)$ and $q$ in place of the original ones $-\nabla V(x)$ and $\Gamma= \sqrt{2D} I$ (compare with (\ref{Eq:HJBequationSmoothSubsolution})). In particular, the second moment of an estimator with change of measure based on the control $v_{t}$ by (\ref{Eq:OptimalControlRough}) will satisfy (\ref{Eq:GoalRegime1Subsolution}); this is Theorem 4.1 in \cite{DupuisSpiliopoulosWang}.

Thus, compared to the case without
multiscale features, one needs to consider the extra factor $\left( I+ \partial \chi(y)/\partial y\right)$, that can be thought as the appropriate {\em weight} function, to achieve asymptotic optimality. In the absence of multiple scales, i.e., when $Q=0$, we have $\chi=0$ and we recover the
case studied in Section \ref{S:EffectRestPoints}. The numerical simulation studies of \cite{DupuisSpiliopoulosWang,DupuisSpiliopoulosWang2} verify the need for accounting for the local environment via the weights $\left( I+ \partial \chi(y)/\partial y\right)$ in the change of simulation measure.

Before illustrating the performance of this importance sampling scheme in a simulation study, let us demonstrate theoretically the necessity to include the cell problem information in the design of the change of measure. For simplicity purposes, let us restrict attention to dimension one. As we have seen before, the effective diffusion coefficient is given by
\[
q=2D \int_{\mathbb{T}}\left(  1+\frac{\partial\chi}{\partial y}\right) ^{2}\mu(dy)
\]

In this case, the optimal change of measure is based on the control
\[
\bar{u}(t,x,y)=-\sqrt{2D}\left( 1+ \partial \chi(y)/\partial y\right) \partial_{x} \bar{U}(t,x).
 \]
 So, let us assume that one is using instead the change of measure, based on the control dictated by the averaged dynamics. Namely, let us assume that the control in question is
$\hat{u}(t,x)=-\sqrt{q}\partial_{x}\bar{U}(t,x)$.

A verification theorem, see \cite{DupuisSpiliopoulosWang} for details,  shows that one would need a statement of the form
\begin{equation}
"\mathbb{E}\int_{t}^{T}\left[  \sqrt{2D}\left(  1+\frac{\partial\chi}{\partial y}\left( \frac
{X^{\epsilon,\delta}_{s}}{\delta}\right)\right)  -\sqrt{q}\right]
ds\rightarrow 0"
\end{equation}

By averaging principle, this is true if
\begin{equation}
\sqrt{q}=\int_{\mathbb{T}}\sqrt{2D}\left(  1+\frac{\partial\chi(y)}{\partial y}\right) \mu(dy).
\end{equation}

However, this is impossible, since
\[
\left(\int\left(  1+\frac{\partial\chi(y)}{\partial y}\right)  \mu(dy)\right)^{2}\neq \int\left(  1+\frac{\partial\chi(y)}{\partial y}\right)^{2} \mu(dy).
\]

This last property explains mathematically why, the local information, as quantified via the cell problem, needs to be taken into account in the design of importance sampling. In Section \ref{SS:SimulationStudyPeriodicCase}, we will also see numerical evidence of this issue.
\subsection{A simulation study}\label{SS:SimulationStudyPeriodicCase}

Let us demonstrate the performance of the importance sampling scheme in a simple simulation study. Consider the one well potential function with diffusion coefficient $D=1$,
\begin{equation}
V(x)=\frac{1}{2}x^2, \hspace{0.5cm}
Q(y)=\cos(y)+\sin(y)
\end{equation}
%\begin{center}
%\includegraphics[height=4cm,width=3cm, angle=90]{Figure1color2}
%\end{center}

Assume that we want to estimate $\theta(\epsilon,\delta)=\mathbb{E}\left[e^{-\frac{1}{\epsilon}h(X^{\epsilon,\delta}_{1})}\right]$, where
%\begin{equation}
 $h(x)=\left(|x|-1\right)^{2}.$
%\end{equation}

It is easy to see that we are dealing with a rare event here, as the function $h(x)$ is minimized at $|x|=1$. Let us compare the following three different estimators
\begin{eqnarray}
\hat{\theta}_{0}(\epsilon,\delta)&=&\frac{1}{K}\sum_{j=1}^{K}\left[e^{-\frac{1}{\epsilon}h(X^{\epsilon,\delta}_{1}(j))}\right] \textrm{ - standard Monte Carlo }\nonumber\\
\hat{\theta}_{1}(\epsilon,\delta)&=&\frac{1}{K}\sum_{j=1}^{K}\left[e^{-\frac{1}{\epsilon}h(\bar{X}^{\epsilon,\delta, \bar{u}}_{1}(j))}Z^{\bar{u}}_{j}\right] \textrm{ - optimal}\nonumber\\
\hat{\theta}_{2}(\epsilon,\delta)&=&\frac{1}{K}\sum_{j=1}^{K}\left[e^{-\frac{1}{\epsilon}h(\bar{X}^{\epsilon,\delta, \hat{u}}_{1}(j))}Z^{\hat{u}}_{j}\right] \textrm{ - ignores local information}\nonumber
\end{eqnarray}

where we have defined the controls
\begin{itemize}
\item{$\bar{u}(t,x,y)=-\sqrt{2}\left( 1+ \partial \chi(y)/\partial y\right)G_{x}(t,x)$--asymptotically optimal.}
\item{$\hat{u}(t,x)=-\sqrt{q}G_{x}(t,x)$--based only on the homogenized system.}
\end{itemize}
and the likelihood ratio is $Z^{u}_{j}=\frac{dP}{d\bar{P}}(\bar{X}^{\epsilon,\delta, u}_{1}(j))$. Notice that in this case, we can compute
\[
1+\frac{\partial\chi(y)}{\partial y}=e^{Q(y)}/\int_{\mathbb{T}} e^{Q(y)}dy,
\]
which justifies the interpretation of the term $1+\frac{\partial\chi(y)}{\partial y}$ as the proper weight term needed that takes into account the local information.

In Table \ref{T:Table3}, we see simulation studies based on $N=10^{7}$ simulation trajectories each, for the estimation of $\theta(\epsilon,\delta)$ using the three different estimators. The measure of comparison is chosen to be the relative error per sample, defined to be
\begin{equation*}
\hat{\rho}_{i}(\epsilon,\delta)\doteq\sqrt{N}\frac{\sqrt{\textrm{Var}(\hat{\theta}_{i}(\epsilon,\delta))}}{\hat{\theta}_{1}(\epsilon,\delta)}.
\end{equation*}

\begin{table}[th]
\begin{center}
\begin{tabular}
[c]{|c|c|c|c|c|c|c|c|c|}\hline
No. & $\epsilon$ & $\delta$ & $\epsilon/\delta$ & $\hat{\theta}_{1}(\epsilon,\delta)$
& $\hat{\rho}_0(\epsilon,\delta)$ & $\hat{\rho}_1(\epsilon,\delta)$ & $\hat{\rho
}_2(\epsilon,\delta)$ \\\hline
$1$ & $0.25$ & $0.1$ & $2.5$ & $2.25e-01$ & $1$ & $6$ & $20$  \\\hline
$2$ & $0.125$ & $0.04$ & $3.125$ & $3.65e-02$ & $3$ &
$6$ & $5$  \\\hline
$3$ & $0.0625$ & $0.015625$ & $4$ & $8.75e-04$ & $34$ &
$4$ & $13$  \\\hline
$4$ & $0.03125$ & $0.007$ & $4.46$ & $6.87e-07$ & $141$ &
$3$ & $105$ \\\hline
$5$ & $0.025$ & $0.004$ & $6.25$ & $1.61e-08$ & $217$ &
$2$ & $97$ \\\hline
$6$ & $0.02$ & $0.002$ & $10$ & $1.99e-10$ & $1294$ & $1$
& $157$  \\\hline
$7$ & $0.015$ & $0.0013$ & $11.54$ & $1.37e-13$ & $800$ &
$1$ & $588$ \\\hline
\end{tabular}
\end{center}\caption{Comparing different importance sampling estimators }\label{T:Table3}
\end{table}

It is clear, that the importance sampling scheme based on the asymptotically optimal change of measure $\bar{u}(t,x,y)$ outperforms the standard Monte Carlo estimator in which no change of measure is being done. It also outperforms,  the estimator based solely on the homogenized system, which ignores the local information characterized by solution to the cell  problem $\chi(y)$.

\section{Importance sampling for multiscale diffusions in random environments}\label{S:RareEventRandomEnvironments}

Let $0<\epsilon,\delta\ll 1$ and consider the process $\left(X^{\epsilon}, Y^{\epsilon}\right)=\left\{\left(X^{\epsilon}_{t}, Y^{\epsilon}_{t}\right), t\in[0,T]\right\}$ taking values in the space $\mathbb{R}^{m}\times\mathbb{R}^{d-m}$ that satisfies the system of SDEs

\begin{eqnarray}
dX^{\epsilon}_{t}&=&\left[  \frac{\epsilon}{\delta}b\left(Y^{\epsilon}_{t},\gamma\right)+c\left(  X^{\epsilon}_{t}%
,Y^{\epsilon}_{t},\gamma\right)\right]   dt+\sqrt{\epsilon}%
\sigma\left(  X^{\epsilon}_{t},Y^{\epsilon}_{t},\gamma\right)
dW_{t},\nonumber\\
dY^{\epsilon}_{t}&=&\frac{1}{\delta}\left[  \frac{\epsilon}{\delta}f\left(Y^{\epsilon}_{t},\gamma\right)  +g\left(  X^{\epsilon}_{t}%
,Y^{\epsilon}_{t},\gamma\right)\right] dt+\frac{\sqrt{\epsilon}}{\delta}\left[
\tau_{1}\left(  Y^{\epsilon}_{t},\gamma\right)
dW_{t}+\tau_{2}\left(Y^{\epsilon}_{t},\gamma\right)dB_{t}\right], \label{Eq:Main}\\
X^{\epsilon}_{0}&=&x_{0},\hspace{0.2cm}Y^{\epsilon}_{0}=y_{0}\nonumber
\end{eqnarray}

We assume non-degeneracy of the diffusion coefficients as well $\mathcal{C}^{1}$ smoothness and boundedness of the drift and diffusion coefficients. Moreover, we assume that $\delta=\delta(\epsilon)\downarrow0$ such that $\epsilon/\delta\uparrow\infty$ as $\epsilon\downarrow0$.  $(W_{t}, B_{t})$ is a $2\kappa-$dimensional standard Wiener process. We assume that for each fixed $x\in\mathbb{R}^{m}$,  $b(\cdot,\gamma), c(x,\cdot,\gamma),\sigma(x,\cdot,\gamma),f(\cdot,\gamma)$, $g(x,\cdot,\gamma), \tau_{1}(\cdot,\gamma)$ and $\tau_{2}(\cdot,\gamma)$ are stationary and ergodic random fields in an appropriate probability space $\left(\Gamma,\mathcal{G},\nu\right)$ with $\gamma\in\Gamma$.

\begin{example}\label{Ex:RandomLangevin}
Notice that if we choose $b(y,\gamma)=f(y,\gamma)=-\nabla_{y}Q(y,\gamma)$  for a periodic function $Q(\cdot)$, $c(x,y,\gamma)=-\nabla_{x}V(x)$,  $\sigma(x,y,\gamma)=\tau_{1}(y,\gamma)=\sqrt{2D}$ and $\tau_{2}(y,\gamma)=0$, and set $y_{0}=x_{0}/\delta$, we then get the Langevin equation (\ref{Eq:LangevinEquation2}). In particular, if we make these choices, then we simply have   $Y^{\epsilon}_{t}=X^{\epsilon}_{t}/\delta$ and the model can be interpreted as  diffusion in the rough potential $\epsilon Q(x/\delta,\gamma)+ V(x)$, where the roughness is dictated by  $Q$. In general, $Q$ may not be modelled as a periodic function. One may model $Q$ as a random field; see the simulation study in Subsection \ref{SS:SimulationStudyRandomCase}.
\end{example}

\subsection{Description of the random environment}\label{SS:DectiptionRandomEnvironment}
The large deviations and importance sampling results for (\ref{Eq:Main}), see \cite{Spiliopoulos2013c,Spiliopoulos2015}, are true under certain assumptions on the random medium that we recall here for convenience.  We assume that there is a group of measure preserving transformations $\{\tau_{y}, y\in\mathbb{R}^{d-m}\}$ acting ergodically on $\Gamma$ that is defined as follows.
\begin{definition}
\label{Def:medium}
\begin{enumerate}
\item {$\tau_{y}$ preserves the measure, namely $\forall y\in\mathbb{R}^{d-m}$
and $\forall A\in\mathcal{G}$ we have $\nu(\tau_{y}A)=\nu(A)$.}

\item {The action of $\{\tau_{y}: y\in\mathbb{R}^{d-m}\}$ is ergodic, that is if
$A=\tau_{y}A$ for every $y\in\mathbb{R}^{d}$ then $\nu(A)=0$ or $1$.}

\item { For every measurable function $f$ on $\left(  \Gamma, \mathcal{G},
\nu\right)  $, the function $(y,\gamma)\mapsto f(\tau_{y}\gamma)$ is
measurable on $\left(  \mathbb{R}^{d-m}\times\Gamma, \mathbb{B}(\mathbb{R}%
^{d-m})\otimes\mathcal{G}\right)  $.}
\end{enumerate}
\end{definition}

Let $\tilde{\phi}$ be a square integrable function in $\Gamma$ and define the operator $T_{y}\tilde{\phi}(\gamma)=\tilde{\phi}(\tau_{y}\gamma)$. The operator $T_{y}\cdot$  is a strongly continuous group of unitary maps in
$L^{2}(\Gamma)$, see \cite{Olla1994}. Denote by $D_{i}$ the infinitesimal generator  of $T_{y}$ in the direction
$i$, which is a closed and densely defined generator, see \cite{Olla1994}.

In order to guarantee that the involved functions are ergodic and stationary random fields on $\mathbb{R}^{d-m}$,
for $\tilde{\phi}\in L^{2}(\Gamma)$, let us define the operator $\phi(y,\gamma)=\tilde{\phi}(\tau_{y}\gamma)$.  Similarly, for a measurable function $\tilde{\phi}:\mathbb{R}^{m}\times\Gamma\mapsto\mathbb{R}^{m}$
we consider the (locally) stationary random field $(x,y) \mapsto \tilde{\phi}(x,\tau_{y}\gamma)=\phi(x,y,\gamma)$. Then, it is guaranteed that $\phi(y,\gamma)$ (respectively $\phi(x,y,\gamma)$) is a stationary (respectively locally stationary) ergodic random field.

The coefficients, $b,c,\sigma,f,g,\tau_{1},\tau_{2}$ of (\ref{Eq:Main}) are defined through this procedure and therefore are guaranteed to be ergodic and stationary random fields. For example in the case of the $c$ drift term,  we start with an $L^{2}(\Gamma)$ function $\tilde{c}(x,\gamma)$ and we define the corresponding coefficients  via the relation $c(x,y,\gamma)=\tilde{c}(x,\tau_{y}\gamma)$.

For every $\gamma\in\Gamma$, let us the operator
\[
\mathcal{L}^{\gamma}=f(y,\gamma)\nabla_{y}\cdot+\text{\emph{tr}}\left[\left(
\tau_{1}(y,\gamma)\tau^{T}_{1}(y,\gamma)+\tau_{2}(y,\gamma)\tau^{T}_{2}(y,\gamma)\right)\nabla_{y}\nabla_{y}\cdot\right] \label{Def:OperatorFastProcess}
\]
which is the infinitesimal generator of a Markov process, say $Y_{t,\gamma}$. Using the Markov process $Y_{t,\gamma}$, we can define the so-called environment process, see \cite{KosyginaRezakhanlouVaradhan,PapanicolaouVaradhan1982,Osada1983,Olla1994}, denoted by $\gamma_{t}$. The environment process $\gamma_{t}$ has continuous
transition probability densities with respect to the $d$-dimensional Lebesgue
measure, see \cite{Olla1994}, and is defined by the equations
\begin{align*}
\gamma_{t}  &  =\tau_{Y_{t,\gamma}}\gamma\label{Eq:EnvironmentProcess}\\
\gamma_{0}  &  =\tau_{y_{0}}\gamma\nonumber
\end{align*}

The infinitesimal generator of the Markov process $\gamma_{t}$ is given by
\[
\tilde{L}=\tilde{f}(\gamma)D\cdot+\text{\emph{tr}}\left[ \left( \tilde{\tau
}_{1}(\gamma)\tilde{\tau}_{1}^{T}(\gamma)+\tilde{\tau
}_{2}(\gamma)\tilde{\tau}_{2}^{T}(\gamma)\right)D^{2}\cdot\right].
\]

In order to simplify the presentation, let us assume that the operator  $\tilde{L}$ is in divergence form. In particular, let us set $\tilde{f}(\gamma)=-D
\tilde{Q}(\gamma)$ and $\tilde{\tau}_{1}(\gamma)=\sqrt{2D}\theta=\text{constant}$ and $\tilde{\tau}_{2}(\gamma)=\sqrt{2D}\sqrt{1-\theta^{2}}=\text{constant}$.

Then, we can write the unique ergodic invariant measure for the environment process $\{\gamma_{t}\}_{t\geq 0}$ in closed form; see \cite{Olla1994, Spiliopoulos2015} for more general case which is not necessarily restricted to the gradient case. Denote by $\mathbb{E}^{\nu}$  the expectation
operator with respect to the measure $\nu$. Then , the  measure $\pi(d\gamma)$ defined  on $(\Gamma,\mathcal{G})$ by%
\[
\pi(d\gamma)\doteq\frac{\tilde{m}(\gamma)}{\mathbb{E}^{\nu}\tilde{m}(\cdot
)}\nu(d\gamma), \text{ with }\tilde{m}(\gamma)=\exp[-\tilde{Q}(\gamma)/D].
\]
is the unique ergodic invariant measure for the environment process $\{\gamma_{t}\}_{t\geq0}$.

Next, we need to define the equivalent to the cell problem in the case of periodic coefficients, also known as the macroscopic problem in the homogenization theory. To do so, we first define $\mathcal{H}^{1}=\mathcal{H}^{1}(\nu)$ to be the Hilbert space  equipped with the
inner product
\[
(\tilde{f},\tilde{g})_{1}=\sum_{i=1}^{d}(D_{i}\tilde{f},D_{i}\tilde{g}).
\]

 Let us consider $\rho>0$ and consider the following problem on $\Gamma$
\begin{equation}
\rho\tilde{\chi}_{\rho}-\tilde{L}\tilde{\chi}_{\rho}=\tilde{b}.
\label{Eq:RandomCellProblem}%
\end{equation}

Under the condition $\tilde{b}\in L^{2}(\nu)$ with $\left\Vert \tilde{b}\right\Vert_{\mathcal{H}^{-1}}<\infty$,  Lax-Milgram lemma, see \cite{Olla1994, KomorowskiLandimOlla2012}, guarantees that
equation (\ref{Eq:RandomCellProblem}) has a unique weak solution in the abstract Sobolev space
$\mathcal{H}^{1}$ or equivalently in $\mathcal{H}^{1}(\pi)$. At this point, we note that in the periodic case one also considers (\ref{Eq:RandomCellProblem}), but one can then take $\rho=0$ given that $b$ averages to zero when is integrated against the invariant measure $\pi$. However, in the random case, (\ref{Eq:RandomCellProblem}) with $\rho=0$ does not necessarily have a well defined solution (even if $b$ averages to zero when is integrated against the invariant measure $\pi$), see for example \cite{KomorowskiLandimOlla2012}.

In the general random case, we consider the equation with $\rho>0$ and in the end, the homogenization theorem is proven by taking appropriate sequences $\rho=\rho(\epsilon)$ such that $\rho(\epsilon)\downarrow 0$ as $\epsilon\downarrow 0$. Taking $\rho\downarrow 0$ is allowed by the following well known properties of the solution to (\ref{Eq:RandomCellProblem}), (see \cite{Olla1994,Osada1983,PapanicolaouVaradhan1982}),
\begin{enumerate}
\item{There is a constant $K$ that is independent of $\rho$ such that
\[
\rho\mathbb{E}^{\pi}\left[  \tilde{\chi}_{\rho}(\cdot)\right]  ^{2}%
+\mathbb{E}^{\pi}\left[  D\tilde{\chi}_{\rho}(\cdot)\right]  ^{2}\leq K
\]}
\item{$\tilde{\chi}_{\rho}$ has an $\mathcal{H}^{1}$
strong limit, i.e., there exists a $\tilde{\chi}_{0}\in\mathcal{H}^{1}(\pi)$
such that
\[
\lim_{\rho\downarrow0}\left\Vert \tilde{\chi}_{\rho}(\cdot)-\tilde{\chi}%
_{0}(\cdot)\right\Vert _{1}=0\quad\text{ and }\quad \lim_{\rho\downarrow0}\rho\mathbb{E}^{\pi}\left[  \tilde{\chi}_{\rho}%
(\cdot)\right]  ^{2}=0.
\]
}
\end{enumerate}

\subsection{Large deviations and importance sampling theory for diffusion in random environments.}\label{SS:LDP_IS_randomEnvironments}
Now that we have defined the random environment and explained its properties, let us review the related large deviations and importance sampling theory from \cite{Spiliopoulos2013c,Spiliopoulos2015}. Set for notational convenience $\tilde{\xi}=D\tilde{\chi}_{0}$.

\begin{theorem}[Theorem 3.5 in \cite{Spiliopoulos2013c}]
\label{T:MainTheorem3} Let $\{\left(X^{\epsilon,\gamma},Y^{\epsilon,\gamma}\right),\epsilon>0\}$ be, for fixed $\gamma\in\Gamma$, the unique strong
solution to (\ref{Eq:Main}). Assume non-degeneracy of the diffusion coefficients as well as $\mathcal{C}^{1}$ smoothness and boundedness of the drift and diffusion coefficients. Consider the regime where $\epsilon,\delta\downarrow 0$ such that $\epsilon/\delta\uparrow\infty$. Then, $\{X^{\epsilon,\gamma},\epsilon>0\}$ converges in probability, almost surely with respect to the random environment  $\gamma\in\Gamma$, as $\epsilon,\delta\downarrow 0$ to the solution of the ODE
\[
d\bar{X}_{t}=r(\bar{X}_{t})dt
\]
and satisfies, almost surely with respect to $\gamma\in\Gamma$, the large deviations principle with rate function
\begin{equation*}
S_{t_{0}T}(\phi)=%
\begin{cases}
\frac{1}{2}\int_{t_{0}}^{T}(\dot{\phi}_{s}-r(\phi_{s}))^{T}q^{-1}(\phi_{s}%
)(\dot{\phi}_{s}-r(\phi_{s}))ds & \text{if }\phi\in\mathcal{AC}%
([t_{0},T]) \text{  and } \phi_{t_{0}}=x_{0}\\
+\infty & \text{otherwise.}%
\end{cases}
\label{Eq:ActionFunctional1}%
\end{equation*}
where
\begin{align*}
r(x)&=\lim_{\rho\downarrow0}\mathbb{E}^{\pi}\left[  \tilde{c}(x,\cdot) +D\tilde{\chi}_{\rho
}(\cdot) \tilde{g}(x,\cdot)\right]  =\mathbb{E}^{\pi}[\tilde{c}(x,\cdot)+\tilde{\xi}%
(\cdot)\tilde{g}(x,\cdot)]\nonumber\\
q(x)&=\lim_{\rho\downarrow0}\mathbb{E}^{\pi}\left[  (\tilde{\sigma}(x,\cdot)+D\tilde{\chi}_{\rho}%
(\cdot)\tilde{\tau}_{1}(\cdot))(\tilde{\sigma}(x,\cdot)+D\tilde{\chi}_{\rho}%
(\cdot)\tilde{\tau}_{1}(\cdot))^{T}+\left(D\tilde{\chi}_{\rho}%
(\cdot)\tilde{\tau}_{2}(\cdot)\right)\left(D\tilde{\chi}_{\rho}%
(\cdot)\tilde{\tau}_{2}(\cdot)\right)^{T}\right]\nonumber\\
&=\mathbb{E}^{\pi}\left[  (\tilde{\sigma}(x,\cdot)+\tilde{\xi}(\cdot)\tilde{\tau}_{1}(\cdot))(\tilde{\sigma}(x,\cdot)+\tilde{\xi}
(\cdot)\tilde{\tau}_{1}(\cdot))^{T}+\left(\tilde{\xi}
(\cdot)\tilde{\tau}_{2}(\cdot)\right)\left(\tilde{\xi}
(\cdot)\tilde{\tau}_{2}(\cdot)\right)^{T}\right]
\end{align*}
and $\rho=\rho(\epsilon)=\frac{\delta^{2}}{\epsilon}$.
\end{theorem}

 Notice that the coefficients $r(x)$ and $q(x)$ are obtained by homogenizing (\ref{Eq:Main}) by taking $\delta\downarrow 0$ with $\epsilon$ fixed. The  form of the action functional can be recognized as the one that would come up when considering large deviations for  the homogenized system. This is also implied by the fact that $\delta$ goes to zero faster than $\epsilon$, since  $\epsilon/\delta\uparrow\infty$.

We also remark here that if $b=0$, then $\chi_{\rho}=0$. In this case  $r(x),q(x)$ take the simplified forms $r(x)=\mathbb{E}^{\pi}[\tilde{c}(x,\cdot)]$ and $q(x)=\mathbb{E}^{\pi}\left[  \tilde{\sigma}(x,\cdot) \tilde{\sigma}(x,\cdot)^{T}\right]$.

Turning now to importance sampling, given controls $u_{1}$ and $u_{2}$ one considers the controlled dynamics under the importance sampling measure $\bar{\mathbb{P}}$
\begin{eqnarray}
d\bar{X}^{\epsilon}_{s}&=&\left[  \frac{\epsilon}{\delta}b\left(\bar{Y}^{\epsilon}_{s},\gamma\right)+c\left(  \bar{X}^{\epsilon}_{s}%
,\bar{Y}^{\epsilon}_{s},\gamma\right)+\sigma\left(  \bar{X}_{s}^{\epsilon},\bar{Y}_{s}^{\epsilon},\gamma\right)  u_{1}(s)\right]   dt+\sqrt{\epsilon}%
\sigma\left(  \bar{X}^{\epsilon}_{s},\bar{Y}^{\epsilon}_{s},\gamma\right)
d\bar{W}_{s}, \nonumber\\
d\bar{Y}^{\epsilon}_{s}&=&\frac{1}{\delta}\left[  \frac{\epsilon}{\delta}f\left(\bar{Y}^{\epsilon}_{s},\gamma\right)  +g\left(  \bar{X}^{\epsilon}_{s}
,\bar{Y}^{\epsilon}_{s},\gamma\right)+\tau_{1}\left(\bar{Y}^{\epsilon}_{s},\gamma\right)u_{1}(s)+
\tau_{2}\left(\bar{Y}^{\epsilon}_{s},\gamma\right)u_{2}(s)\right]   dt\nonumber\\
& &\hspace{5cm}+\frac{\sqrt{\epsilon}}{\delta}\left[
\tau_{1}\left(\bar{Y}^{\epsilon}_{s},\gamma\right)
d\bar{W}_{s}+\tau_{2}\left(\bar{Y}^{\epsilon}_{s},\gamma\right)d\bar{B}_{s}\right],\label{Eq:Main2}\\
\bar{X}^{\epsilon}_{t_{0}}&=&x_{0},\hspace{0.2cm}\bar{Y}^{\epsilon}_{t_{0}}=y_{0}\nonumber
\end{eqnarray}

where $(v_{1}(s), v_{2}(s))$ denote the first and second component of the control
\[
u(s,\bar{X}^{\epsilon}_{s},\bar{Y}^{\epsilon}_{s})=(u_{1}(s,\bar{X}^{\epsilon}_{s},\bar{Y}^{\epsilon}_{s}), u_{2}(s,\bar{X}^{\epsilon}_{s},\bar{Y}^{\epsilon}_{s})).
\]

Then, for a given cost function $h(x)$, under $\bar{\mathbb{P}}$
\[
\Delta^{\epsilon,\gamma}(t_{0},x_{0},y_{0})=\exp\left\{  -\frac{1}{\epsilon}h(\bar{X}^{\epsilon
}_{T})\right\}  \frac{d\mathbb{P}}{d\bar{\mathbb{P}}}(\bar{X}^{\epsilon
}, \bar{Y}^{\epsilon}),
\]
is an unbiased estimator for
$\mathbb{E}\left[\exp\left\{  -\frac{1}{\epsilon}h(X^{\epsilon}_{T})\right\}\right]$.

Consider next the Hamiltonian
\[
H(x,p)=\left<r(x), p\right>-\frac{1}{2}\left\Vert q^{1/2}(x)p\right\Vert^{2}
\]
 with $r(x),q(x)$ the coefficients defined in Theorem \ref{T:MainTheorem3} and consider the HJB equation associated to this Hamiltonian, letting $\bar{U}(t,x)$ be a smooth subsolution to it (analogously to Section \ref{S:ReviewIS_LDP} with $r(x)$ and $q(x)$ in place of $-\nabla V(x)$ and $\Gamma$ respectively). Then, the following theorem guarantees at least logarithmic asymptotically good performance.

\begin{theorem}[Theorem 4.1 in \cite{Spiliopoulos2015}]
\label{T:UniformlyLogEfficient} Let $\{\left(X^{\epsilon}_{s}, Y^{\epsilon}_{s} \right),\epsilon>0\}$ be
the solution to (\ref{Eq:Main}) for $s\in[t_{0},T]$ with initial point $(x_{0},y_{0})$ at time $t_{0}$.  Consider a non-negative, bounded and continuous
function $h:\mathbb{R}^{m}\mapsto\mathbb{R}$. Let $\bar{U}(s,x)$ be a subsolution to the associated HJB equation that has continuous derivatives up to
order $1$ in $t$ and order $2$ in $x$, and the first and second
derivatives in $x$ are uniformly bounded. Assume non-degeneracy of the diffusion coefficients as well $\mathcal{C}^{1}$ smoothness and boundedness of the drift and diffusion coefficients. In the general case where $b\neq 0$, consider $\rho>0$ and define the (random) feedback control $u_{\rho}(s,x,y,\gamma)=\left(u_{1,\rho}(s,x,y,\gamma), u_{2,\rho}(s,x,y,\gamma)\right)$ by

 \begin{equation*}
u_{\rho}(s,x,y,\gamma)=\left(-\left(\sigma+D\chi_{\rho} \tau_{1}\right)^{T}(x,y,\gamma)\nabla_{x}\bar{U}(s,x), -\left(D\chi_{\rho}\tau_{2}\right)^{T}(y,\gamma)\nabla_{x}\bar{U}(s,x)\right)\label{Eq:feedback_controlReg1}
\end{equation*}

Then for $\rho=\rho(\epsilon)=\frac{\delta^{2}}{\epsilon}\downarrow 0$ we have that almost surely in $\gamma\in\Gamma$
\begin{equation}
\liminf_{\epsilon\rightarrow0}-\epsilon\ln Q^{\epsilon,\gamma}(t_{0},x_{0},y_{0};u_{\rho}(\cdot))\geq
G(t_{0},x_{0})+\bar{U}(t_{0},x_{0}). \label{Eq:GoalSubsolution}%
\end{equation}

If $b=0$, then set $u(s,x,y,\gamma)=\left(-\sigma^{T}(x,y,\gamma)\nabla_{x}\bar{U}(s,x), 0\right)$ and (\ref{Eq:GoalSubsolution}) holds with $u_{\rho}(\cdot)=u(\cdot)$.
\end{theorem}

%For the sake of completeness, we also mention here that Theorem \ref{T:UniformlyLogEfficient} can be also used to obtain  good asmptotic performance for estimating probabilities. For instance, if we define
%\[
%h(x)=%
%\begin{cases}
%0 & \text{if }x\in A\\
%+\infty & \text{if }x\notin A,
%\end{cases}
%\]
%then Theorem \ref{T:UniformlyLogEfficient} is dealing with probabilities of the form $\mathrm{P}_{t_{0},x_{0},y_{0}}[X^{\epsilon,\gamma}_{T}\in A]$.

\subsection{A simulation study}\label{SS:SimulationStudyRandomCase}
Consider for instance the case of Example \ref{Ex:RandomLangevin}
\begin{equation}
dX^{\epsilon,\delta}_{t}=-\nabla
V^{\epsilon}\left(X^{\epsilon,\delta}_{t},\frac{X^{\epsilon,\delta}_{t}}{\delta}\right)dt+\sqrt{2\epsilon}dW_{t},
\end{equation}
where the potential function
$V^{\epsilon}\left(x,x/\delta\right)=\epsilon Q(x/\delta)+V(x).$
$Q(y)$ is a stationary ergodic random field on a probability space $(\mathcal{X},\mathcal{G},\nu)$. We may consider for instance $V(x)=\frac{1}{2}x^2$ and
\begin{equation*}
Q(y)\textrm{ mean zero Gaussian with } \mathrm{E}^{\nu}\left[Q(x)Q(y)\right]=\exp\left[-\left|x-y\right|^{2}\right]
\end{equation*}

Making the connection with (\ref{Eq:Main}), the fast $Y$ motion essentially is $Y=X/\delta$. Referring to Theorems \ref{T:MainTheorem3} and \ref{T:UniformlyLogEfficient} we have $r(x)=-V^{\prime}(x)/(K\hat{K})$ and $q=2/(K\hat{K})$ where $K=\mathrm{E}^{\nu}[  e^{-Q(z)}], \hat{K}=\mathrm{E}^{\nu}[  e^{Q(z)}]$. Given a classical subsolution $\bar{U}$, one expects that the corresponding change of simulation measure that guarantees at least asymptotic optimality, is based on the control
$\bar{u}(s,x,y,\gamma)=(-\sqrt{2}(1+\partial \chi(y,\gamma)/\partial y)\bar{U}_{x}(s,x),0)$ where one can compute that the weight function is $1+\partial \chi(y,\gamma)/\partial y=e^{Q(y,\gamma)}/\hat{K}$. Note that in contrast to the periodic case, the control $u$ is random in that it implicitly depends on $\gamma\in\Gamma$, via the random field $Q(y,\gamma)$.

Assume that we want to estimate
\begin{equation}
 \theta^{\epsilon,\delta}=\mathrm{P}\left[X^{\epsilon,\delta} \textrm{ hits } 1 \textrm{ before } 0|X^{\epsilon,\delta}_{0}=0.1\right]
\end{equation}

As in Subsection \ref{SS:SimulationStudyPeriodicCase}, we compare the asymptotical optimal change of measure with standard Monte Carlo, which corresponds to no change of measure, and with the importance sampling that corresponds to the change of measure based only on the homogenized problem, which ignores the macroscopic problem.
Based on $10^{7}$ trajectories, we have the following simulation data
\begin{table}[th]
\begin{center}
\begin{tabular}
[c]{|c|c|c|c|c|c|c|c|}\hline
No. & $\epsilon$ & $\delta$ & $\epsilon/\delta$ & $\hat{\theta}_{1}(\epsilon,\delta)$  & $\hat{\rho
}_{0}(\epsilon,\delta)$ & $\hat{\rho}_{1}(\epsilon,\delta)$ & $\hat{\rho}_{2}(\epsilon,\delta
)$\\\hline
$1$ & $0.25$ & $0.1$ & $2.5$  & $1.38e-1$  & $3$ & $0.5$ & $3$\\\hline
$2$ & $0.125$ & $0.04$ & $3.125$ &  $1.31e-2$  & $7$ & $16$ & $8$\\\hline
$3$ & $0.0625$ & $0.018$ & $3.472$  & $6.13e-4$  & $36$ & $18$ & $42$\\\hline
$4$ & $0.05$ & $0.01$ & $5$ &  $2.30e-5$ & $212$ & $28$ & $316$\\\hline
$5$ & $0.04$ & $0.007$ & $5.72$  & $5.93e-6$  & $396$ & $75$ & $332$\\\hline
$6$ & $0.025$ & $0.004$ & $6.25$  & $7.82e-10$  & $-$ & $22$ & $1856$\\\hline
\end{tabular}
\end{center}\caption{Comparing different importance sampling estimators with $x^{-}=0$ (equilibrium), $x_{0}=0.1$ (initial point), $x^{+}=1$ (target).}\label{T:Table4}
\end{table}

It is clear, that the importance sampling scheme based on the asymptotically optimal change of measure $\bar{u}(t,x,y,\gamma)$ outperforms the standard Monte Carlo estimator in which no change of measure is being done. It also outperforms,  the estimator based solely on the homogenized system, which ignores the local information characterized by solution to the macroscopic problem. Of course, this behavior is parallel to the behavior observed in the periodic case of Subsection \ref{SS:SimulationStudyPeriodicCase}. Additional simulation studies can be found in \cite{Spiliopoulos2015}.

In \cite{Spiliopoulos2015}, the interested reader can find further simulation studies in the case of the general model (\ref{Eq:Main}) where one does not necessarily have the $Y$ motion to be $X/\delta$. However, we do point out that the theoretical results of \cite{Spiliopoulos2015} are valid for the system (\ref{Eq:Main}) where the process $(X^{\epsilon},Y^{\epsilon})$ has initial point $(x_{0},y_{0})$ and both $x_{0}$ and $y_{0}$ are of order one as $\delta\downarrow 0$. This is not exactly the same to the case where $Y=X/\delta$, as then $y_{0}=x_{0}/\delta$, which is no longer of order one as $\delta\downarrow 0$. But, simulation studies, as the one presented in Table \ref{T:Table4}, indicate that the theoretical results should be also valid for the $Y=X/\delta$ case.

\section{Importance sampling for metastable multiscale processes and further challenges}\label{S:IS_metastabilityMultiscale}
In Section \ref{S:EffectRestPoints} we elaborated on the effects of rest points and metastable dynamics on importance sampling schemes. The end conclusion was that extra care is needed when stable or unstable equilibrium points are in the domain of interest. In this case, asymptotic optimality is not enough in that asymptotically optimal schemes may not perform well in practice unless one goes to really small values of $\epsilon$, in which case the events may be too rare to be of any practical interest. Then, in Section \ref{S:IS_MultiscaleLangevin} and \ref{S:RareEventRandomEnvironments} we summarized the issues that come up in the design of asymptotically efficient importance sampling schemes when the dynamics have multiple scales.

In \cite{DupuisSpiliopoulosZhou2013,DupuisSpiliopoulos2014a} we have systematically addressed the effects of rest points onto the design of importance sampling schemes and have identified what the main issues are. In \cite{DupuisSpiliopoulosZhou2013}, we have suggested a potential provably appropriate remedy to the issue, by constructions as the ones mentioned in Section \ref{S:EffectRestPoints}. The subsolution constructed there effectively yields a very good approximation to the zero variance change of measure. %The same construction with a good pre-asymptotic performance to the same degree is more difficult to be carried over in higher dimensions. This point is being further discussed in \cite{DupuisSpiliopoulos2014a}, where such constructions are demonstrated to work  reasonably well in higher dimensions.
Even though, the constructions in \cite{DupuisSpiliopoulosZhou2013,DupuisSpiliopoulos2014a} work provably well pre-asymptotically and asymptotically and do not degrade as parameters such as the time horizon $T$ getting large, the performance in higher dimensions can be worse than the corresponding  performance in the lower-dimensional cases. While this is expected to be the case as the dimension gets larger, due to further approximations and simplifications that need to be made, there is a clear room for improvement here. This is part of ongoing work of the author and we refer the interested reader to \cite{SalinsSpiliopoulos2016} for some results in the infinitely dimensional small noise SPDE case.

Moreover, it is clear that the constructions of Sections \ref{S:IS_MultiscaleLangevin} and \ref{S:RareEventRandomEnvironments} guarantee only asymptotic optimality. If in addition to multiscale dynamics one has to also face metastability, then, as it was seen in Section \ref{S:EffectRestPoints}, theoretical asymptotic optimality is not sufficient for good numerical performance. One can of course combine the results of Section \ref{S:EffectRestPoints} with those of Sections \ref{S:IS_MultiscaleLangevin} and \ref{S:RareEventRandomEnvironments}. To be more precise, one can combine the results of \cite{DupuisSpiliopoulosZhou2013,DupuisSpiliopoulos2014a} with those of \cite{DupuisSpiliopoulosWang,Spiliopoulos2015}. In practice, one can just use the changes of measure as indicated in \cite{DupuisSpiliopoulosWang,Spiliopoulos2015} that guarantee asymptotic optimality, but construct the subsolution $\bar{U}(t,x)$ as indicated in \cite{DupuisSpiliopoulosZhou2013,DupuisSpiliopoulos2014a}. We plan to address this issue in more detail in a future work.

%\bibliographystyle{plain}
%\bibliography{main}

%%%%%%%%%%%%%%%%%%%%%%%%%%%%%%%%%%%%%%%%%%%%%%%%%%%%%%%%%%%%%
%\begin{thebibliography}{99}
%%%%%%%%%%%%%%%%%%%%%%%%%%%%%%%%%%%%%%%%%%%%%%%%%%%%%%%%%%%%%%%%%%
%

%\end{thebibliography}

\end{document}